\documentclass[journal]{IEEEtran}

\usepackage{amsmath,amsfonts}
\usepackage{algorithm}
\usepackage{algorithmic}
\usepackage{array}
\usepackage[caption=false,font=normalsize,labelfont=sf,textfont=sf]{subfig}
\usepackage{textcomp}
\usepackage{stfloats}
\usepackage{url}
\usepackage{verbatim}
\usepackage{graphicx}
\usepackage{cite}
\usepackage{bm}
\usepackage{amsthm}
\newtheorem{theorem}{Theorem}
\newtheorem{lemma}{Lemma}
\newtheorem{proposition}{Proposition}
\newtheorem{corollary}{Corollary}

\newtheorem{remark}{Remark}

\usepackage{balance}
\usepackage{etoolbox}

\makeatletter
\def\@IEEEBIOskipN{1\baselineskip}
\expandafter\patchcmd\csname\string\IEEEbiography\endcsname
  {\vskip \@IEEEBIOskipN plus 1fil minus 0\baselineskip}
  {\vskip \@IEEEBIOskipN\relax}
  {\typeout{** BIOGAP PATCH OK **}}
  {\typeout{** BIOGAP PATCH FAILED **}}
\makeatother

\begin{document}

\title{A Bilevel Neurodynamic Approach for Distributed SCAD-Regularized Sparse Least Squares}

\author{Juntao Li, Qihang Guo, and Deyuan Meng%
\thanks{J. Li and Q. H. Guo are with the School of Mathematics and Statistics, Henan Normal University, Xinxiang 453007, P. R. China (e-mail: juntaolimail@126.com).}%
\thanks{D. Meng is with the Seventh Research Division, School of Automation Science and Electrical Engineering, Beihang University (BUAA), Beijing 100191, P. R. China (e-mail: dymeng@buaa.edu.cn).}%
\thanks{This work was supported in part by the National Natural Science Foundation of China under Grants U2333215 and 62273018.}%
\thanks{Corresponding author: Deyuan Meng.}%
}

\maketitle

\begin{abstract}
This paper proposes a bilevel particle swarm optimization (PSO)-neurodynamic method for solving the distributed sparse least-squares problem regularized by the smoothly clipped absolute deviation (SCAD) penalty to global optimality. The SCAD penalty is nonsmooth and nonconvex, so existing distributed neurodynamic methods have not been extended to it directly. Through two-time-scale singular perturbation analysis, a rigorous full-trajectory convergence theory with explicit convergence rates is established for the inner collaborative multi-agent neurodynamic system without requiring the second-order sufficient condition. A theory-guided inner-outer collaborative optimization framework is further developed, in which a red-light green-light sensitivity mechanism is embedded in the outer PSO search to enhance global search while concentrating computational effort on the SCAD shape parameters, thereby improving recovery accuracy and computational efficiency.
Simulation results on four sparse-recovery scenarios demonstrate that the proposed method accurately recovers the true sparsity pattern and achieves faster convergence and lower relative error than competing neurodynamic methods.
\end{abstract}


\section{Introduction}
\label{sec:introduction}

Sparse least-squares recovery is a central
problem in signal processing and system
identification, where the objective is to
recover a sparse signal from underdetermined
linear
measurements~\cite{donoho2006compressed,
candes2006robust}.  The convex \(\ell_1\)
relaxation~\cite{tibshirani1996regression}
is the standard computationally tractable
approach.  Centralized algorithms such as
gradient projection, the iterative
shrinkage-thresholding algorithm
(ISTA)~\cite{daubechies2004iterative} and
its fast variant
(FISTA)~\cite{beck2009fast} and the
alternating direction method of multipliers
(ADMM)~\cite{boyd2011distributed} have been
developed for the
\(\ell_1\)-regularized problem, but these
methods assume that all observations are
available at a single processing node.  In
large-scale networked systems, relaying all
measurements to a central fusion center
incurs prohibitive communication overhead
and potential privacy risks.  This practical
limitation has motivated distributed
solution paradigms such as consensus-based
ADMM~\cite{boyd2011distributed} and
distributed proximal gradient methods for
weakly convex
regularization~\cite{komuro2022distributed},
which enable individual nodes to
collaboratively reconstruct sparse signals
using only local measurements and
neighbor-to-neighbor communication.

The convex \(\ell_1\) relaxation, however,
introduces bias on large coefficients.  To
mitigate this limitation, nonconvex
regularizers such as the minimax concave
penalty~\cite{zhang2010mcp} have been
introduced, among which the smoothly clipped
absolute deviation penalty
function~\cite{fan2001variable}, hereafter
abbreviated as the SCAD penalty function, is
a prominent example.  The SCAD penalty
function achieves higher-fidelity
recovery, but replacing
the convex \(\ell_1\) norm with it
introduces two coupled computational
challenges, namely nonsmoothness and
nonconvexity.  The SCAD penalty function is
nondifferentiable at the origin and at the
transition points determined by its shape
parameters, so its gradient mapping lacks
global Lipschitz continuity.  It is concave
on an interval between two shape
parameters, which creates a
landscape populated with spurious local minima
and causes severe initialization
sensitivity~\cite{rockafellar1997variational,
xue2012nonconcave}.

To address the nonsmoothness challenge,
smoothing techniques approximate the
penalty with a differentiable surrogate.
Nesterov's smoothing
technique~\cite{nesterov2005smooth} provides
the foundational framework by approximating
nonsmooth functions via convolution with a
smooth kernel.  This idea has been extended
to distributed settings, where the
decentralized smoothing ADMM
framework~\cite{mirzaeifard2025decentralized}
applies a Nesterov-type smoothing kernel to
approximate the penalty and solves the
resulting problem distributively.
Alternative smoothing approaches include the
Moreau envelope for weakly convex
regularizers in distributed sparse
optimization~\cite{komuro2022distributed}
and duality-based smoothing
methods~\cite{antil2026duality}.  These
techniques convert the nondifferentiable
penalty into a smooth one, but the uniform
kernel injects a nonzero derivative in the
large-coefficient region where the original
SCAD penalty function is constant,
destroying its piecewise structure
and introducing bias on large coefficients.
Moreover, these techniques are iterative and
offer no continuous-time transient analysis,
while the nonconvexity of the landscape
remains unresolved, leaving the second
challenge to a different paradigm.

For the nonconvexity, continuous-time
neurodynamic optimization offers this
alternative by evolving the
decision variables as a continuous-time
dynamical system toward minimizers of the
nonconvex objective.  Originating with the
Hopfield networks~\cite{hopfield1982neural}
and the Hopfield-Tank
networks~\cite{hopfield1985neural}, the
classical projection neural network
framework~\cite{xia2004projection,
hu2007gpnn} formalized
this paradigm for constrained convex
optimization, and it has since been
advanced by trust-region and
noise-suppression projection network
designs~\cite{huang2025trustregion,
yang2025noisesuppression}.  As problems scale to
multi-agent networks, distributed methods
have emerged, starting from a collective
neurodynamic approach to distributed
constrained optimization~\cite{liu2017collective}
and a collaborative neurodynamic framework
that coordinates multiple recurrent
networks~\cite{yang2018collaborative}.
These efforts have led to the projection
neural network with an augmented Lagrangian
~\cite{xia2023collaborative},
the two-timescale projection neural network
~\cite{huang2024twotimescale},
the penalty neurodynamic approach
~\cite{jia2024collective},
and a momentum-based
scheme~\cite{xia2024momentum}.  For
problems with multiple local minima,
collaborative neurodynamic optimization
(CNO)
frameworks~\cite{che2019collaborative,
che2021twotimescale, wang2024supervised}
embed a neurodynamic
inner solver within a metaheuristic outer
search.  Recent developments have extended
these paradigms toward sparse recovery and
nonsmooth composite optimization, including
distributed projection networks for
block-decomposed basis pursuit~\cite{zhao2025distributed},
proximal gradient neurodynamic networks for
sparse optimization~\cite{xu2025fixedtime},
collaborative neurodynamic formulations for
compressed sensing~\cite{li2026compressed},
and proximal neural networks for nonsmooth
composite problems~\cite{luan2026antidisturbance}.
These methods are designed for
objectives with globally Lipschitz gradients,
an assumption that has proven effective
across a broad class of nonconvex problems;
the nonsmoothness of the SCAD penalty,
however, lies outside this setting.  When the
penalty is smoothed to restore
differentiability, their Lyapunov analyses
typically invoke the second-order sufficient
condition (SOSC) that the Hessian is
positive definite on the active subspace at
the equilibrium.  These analyses characterize
the asymptotic phase as time approaches
infinity within an equilibrium neighborhood,
and the boundary-layer transient by which
trajectories enter this neighborhood remains
to be further characterized.  Since this
transient governs the parameter-to-trajectory
mapping required by an outer search, CNO
frameworks have so far selected
hyperparameter dimensions by heuristic rules
rather than through the inner convergence
theory.

Although the neurodynamic methods reviewed
above offer a powerful framework for
nonconvex optimization, they have not yet
been applied to the SCAD-regularized
problem, whose penalty function is
nonsmooth.  To the best of our
knowledge, addressing the
nonsmoothness and the nonconvexity of the
SCAD penalty function simultaneously remains
an open challenge for distributed
neurodynamic methods, even as distributed
nonsmooth nonconvex optimization continues to
attract attention~\cite{liu2026distributed}.  A
structure-preserving mollification of the
SCAD penalty is therefore introduced to
yield a smooth surrogate gradient, by which
the second-order sufficient condition and
local strong convexity are obtained near a
strict local minimum.  This paper develops
a bilevel distributed particle swarm
optimization (PSO)-neurodynamic method
around this mollification, in which the
inner neurodynamic system achieves
full-trajectory convergence without the
SOSC and the outer PSO search, guided by
the inner convergence rates, mitigates
the spurious local minima caused by
nonconvexity.  The main
contributions are summarized as follows.
\begin{enumerate}
    \item A bilevel PSO-neurodynamic method is proposed
    for solving the distributed
    SCAD-regularized sparse least-squares
    problem to global optimality.  The
    inner level is a collaborative
    multi-agent neurodynamic system that
    uses only neighbor-to-neighbor communication to
    achieve local convergence, whereas the
    outer level is a particle swarm search
    layer that mitigates the spurious
    local minima problem and
    initialization sensitivity caused by
    the nonconvexity of SCAD. The proposed approach does not rely on the SOSC invoked by existing nonconvex neurodynamic methods \cite{xia2023collaborative, huang2024twotimescale, jia2024collective, xia2024momentum, che2019collaborative}, thereby relaxing the applicability conditions.

    \item A rigorous full-trajectory convergence theory is established for the inner multi-agent neurodynamic system via two-time-scale singular perturbation analysis. The system trajectory is shown to converge exponentially to the slow manifold at rate \(\alpha\) during the boundary-layer transient \([0,T]\), and to an \(O(\sqrt{\epsilon_1})\)-neighborhood of the equilibrium set at rate \(\zeta/(2k_1)\) during the slow phase \([T,\infty)\).
         Existing distributed neurodynamic analyses~\cite{xia2023collaborative,huang2024twotimescale,jia2024collective,xia2024momentum} establish convergence within a neighborhood of the equilibrium point; the proposed theory goes further by characterizing the complete trajectory from the initial state to the equilibrium neighborhood and by providing an explicit convergence rate, thereby offering a quantitative foundation for outer global optimization design.

    \item A theory-guided inner--outer collaborative optimization framework is proposed, in which a red-light green-light (RLGL) sensitivity mechanism is embedded in the outer particle swarm search to enhance its global search capability. The RLGL mechanism uses the inner convergence rate \(\zeta/(2k_1)\) as a sensitivity surrogate to freeze dimensions with negligible sensitivity in the outer search, while the theory-derived constraint \(\epsilon_1<\epsilon_1^*\) is incorporated into the outer search to ensure the validity of the inner singular perturbation analysis and the stability of the inner neurodynamic system. In contrast to frameworks that connect the two layers through hyperparameter tuning~\cite{xia2023collaborative,huang2024twotimescale,jia2024collective,xia2024momentum,che2019collaborative}, the proposed framework preserves sparsity selection and focuses the computational effort on the SCAD shape parameters, thereby improving both recovery accuracy and computational efficiency.
\end{enumerate}

The remainder is organized as follows.
Section~\ref{sec:problem} formulates the problem.
Section~\ref{sec:main_results} develops the
structure-preserving mollification and the
collaborative neurodynamic system.
Section~\ref{sec:convergence} establishes the
full-process convergence theory.
Section~\ref{sec:swarm} presents the bilevel
PSO-neurodynamic approach with the RLGL
mechanism.
Section~\ref{sec:simulation} presents numerical
experiments, and Section~\ref{sec:conclusion}
concludes.

\section{Problem Description and Preliminaries}
\label{sec:problem}

Sparse least-squares optimization, which seeks
sparse solutions to linear inverse problems,
arises in diverse domains including machine
learning, signal processing, and system
identification.  The problem of finding an
intrinsically sparse vector \(\mathbf{x} \in
\mathbb{R}^n\) from linear measurements is
formally stated as follows. Given a measurement
matrix \(\mathbf{A} \in \mathbb{R}^{m \times n}\),
an observation vector \(\mathbf{b} \in \mathbb{R}^m\),
and additive noise \(\mathbf{e} \in \mathbb{R}^m\),
this paper is devoted to solving the
SCAD-regularized sparse least-squares problem:
\begin{equation}
    \min_{\mathbf{x} \in \mathbb{R}^n} \; \frac{1}{2} \|\mathbf{A} \mathbf{x} - \mathbf{b}\|^2 + P_{\lambda}(\mathbf{x}),
    \label{eq:sparse_lsq}
\end{equation}
where \(P_{\lambda}(\mathbf{x})  = \sum_{j=1}^{n}
p_{\lambda}(x_j)\) denotes the SCAD penalty, the
univariate component \(p_{\lambda}(x_j)\) equals
\(\lambda |x_j|\) if \(|x_j| \le \lambda\),
\(-(|x_j|^2 - 2a\lambda |x_j| +
\lambda^2)/(2(a-1))\) if \(\lambda < |x_j| \le
a\lambda\), and \((a+1)\lambda^2/2\) if \(|x_j| >
a\lambda\), with regularization parameter
\(\lambda > 0\) and shape parameter \(a
> 2\).  The SCAD penalty is both
\emph{nonconvex} (it is concave on
\((\lambda,a\lambda)\)) and
\emph{nonsmooth} (it is nondifferentiable at
\(0,\pm\lambda,\pm a\lambda\)), which makes
\eqref{eq:sparse_lsq} a nonconvex, nonsmooth
optimization problem admitting multiple local
minima and precluding direct use of
gradient-based dynamics.

Before proceeding, we fix the notation
conventions adopted in the rest of the paper. Let
\(\mathbb{R}^n\) denote the \(n\)-dimensional
Euclidean space, and let \(\mathbf{I}_n\) and
\(\mathbf{1}_N\) be the \(n\times n\) identity
matrix and the \(N\)-dimensional all-ones vector,
respectively. For any vector \(\mathbf{x}\),
\(\|\mathbf{x}\|\) denotes its Euclidean norm.
The operators \(\operatorname{col}(\cdot)\) and
\(\operatorname{null}(\cdot)\) represent
column-wise concatenation and matrix null space,
respectively, and \(\otimes\) is the Kronecker
product. For a symmetric matrix,
\(\lambda_{\min}(\cdot)\) and
\(\lambda_{\max}(\cdot)\) return its smallest and
largest eigenvalues.

\section{SCAD Mollification and
Collaborative Neurodynamic System}
\label{sec:main_results}

The original SCAD penalty is nonsmooth at the
origin, preventing the use of standard
gradient-based dynamics. According
to~\cite{rockafellar1997variational,nesterov2005smooth},
We employ the following mollifier
\[
\rho_{\varepsilon}(t)=\varepsilon^{-1}\rho(t/\varepsilon),
\]
where \(\varepsilon>0\) is a parameter to be determined, and \(\rho\) is defined by
\[
\rho(u)=
\begin{cases}
c\,\mathrm{e}^{1/(u^2-1)}, & |u|<1,\\[4pt]
0, & |u|\ge 1,
\end{cases}
\]
with normalization constant
$
c=\left(\int_{-1}^{1}\mathrm{e}^{1/(u^2-1)}\,du\right)^{-1}.
$
It is noted that the admissible range and optimal value of parameter \(\varepsilon\) are subsequently determined by the local-optimality theory (Theorem~\ref{thm:local_optimality}) and the outer PSO search (Section~\ref{sec:swarm}), respectively.
Further, we propose a \emph{structure-preserving}
mollified SCAD penalty that---unlike
Moreau--Yosida regularization---convolves only the
derivative \(p'_\lambda\) rather than the penalty
itself, so that the piecewise affine--quadratic
structure of SCAD is retained up to an
\(\varepsilon\)-neighborhood of the transition
points:
\begin{equation}
    P_{\lambda}^{\varepsilon}(\mathbf{x}) = \sum_{j=1}^{n} p_{\lambda}^{\varepsilon}(x_j),
    \label{eq:P_mollified}
\end{equation}
where \(p_{\lambda}^{\varepsilon}(x_{j}) =
\int_{0}^{|x_{j}|} q_{\varepsilon}(t)\,dt\),
$q_{\varepsilon}(t) = \int_{\mathbb{R}}
\rho_{\varepsilon}(s) \, p'_{\lambda}(t - s)  \,
ds$, $p'_{\lambda}(\cdot)$ is the derivative of
SCAD penalty (the value of $p'(*)$ can be
assigned arbitrarily at the nondifferentiable
points). Replacing SCAD penalty in
(\ref{eq:sparse_lsq}) by
$P_{\lambda}^{\varepsilon}$ in
(\ref{eq:P_mollified}) yields the smooth
SCAD-regularized sparse least-squares problem:
\begin{equation}
    \min_{\mathbf{x} \in \mathbb{R}^n} \; \frac{1}{2} \|\mathbf{A} \mathbf{x} - \mathbf{b}\|^2 + P_{\lambda}^{\varepsilon}(\mathbf{x}).
    \label{eq:mollified_central}
\end{equation}

We first establish the Lipschitz continuity of the mollified SCAD gradient, which provides the foundation for the subsequent convergence analysis.

\begin{lemma}
\label{lem:lip} For a mollification parameter \(\varepsilon>0\) to be determined,
 the gradient $\nabla
P_{\lambda}^{\varepsilon}(\mathbf{x})$ is
globally Lipschitz on $\mathbb{R}$ with the
constant $l_\varepsilon = {\frac{2\lambda
c}{\varepsilon\, e} }$.
\end{lemma}

\begin{proof}
Observe that $\rho(u)$ is even and strictly
decreasing on $(0,1)$.  Hence, $\rho'(u)<0$ on
$(0,1)$. It follows from the definition of
the mollifier that
\begin{equation}\label{A}
\begin{aligned}
\int_{\mathbb{R}} |\rho_\varepsilon'(s)|\,ds
&= \frac{1}{\varepsilon^2}\int_{\mathbb{R}} \bigl|\rho'(s/\varepsilon)\bigr|\,ds
= \frac{1}{\varepsilon}\int_{-1}^{1} |\rho'(u)|\,du  \\
&= \frac{2}{\varepsilon}\int_{0}^{1} \bigl(-\rho'(u)\bigr)\,du
= \frac{2c}{\varepsilon e}.
\end{aligned}
\end{equation}
Observe that $\rho_\varepsilon(t)\in C_c^\infty$ and
$p'_\lambda(t)$ is bounded and measurable. Hence,
differentiating the function $q_{\varepsilon}(t)$
and using~\eqref{A} yields
\[
|q_\varepsilon'(t)|
= \biggl|\int_{\mathbb{R}} \rho_\varepsilon'(s)\,p'_\lambda(t-s)\,ds\biggr|
\le \lambda \int_{\mathbb{R}} |\rho_\varepsilon'(s)|\,ds
=l_\varepsilon .
\]
The mean value theorem implies
$|q_\varepsilon(t_1)-q_\varepsilon(t_2)|\le
l_\varepsilon|t_1-t_2|$ for all $t_1,t_2$.
Observe that \(\nabla P_\lambda^\varepsilon(\mathbf{x}) = \bigl(q_\varepsilon(x_1),\dots,q_\varepsilon(x_n)\bigr)^T\) and each component satisfies Lipschitz continuity. Hence,
$
\|\nabla P_\lambda^\varepsilon(\mathbf{x}) - \nabla P_\lambda^\varepsilon(\mathbf{y})\|_2 \le l_\varepsilon \|\mathbf{x}-\mathbf{y}\|_2, \quad \forall \mathbf{x},\mathbf{y}\in\mathbb{R}^n.
$
\end{proof}

To solve the smooth
problem~\eqref{eq:mollified_central} in a fully
distributed manner, we consider a network of
\(N\) agents communicating over a connected
undirected graph
\(\mathcal{G}=(\mathcal{V},\mathcal{E})\), where
\(\mathcal{V}=\{1,\dots,N\}\) is the vertex set
and
\(\mathcal{E}\subseteq\mathcal{V}\times\mathcal{V}\)
is the edge set. The neighbor set of agent \(i\)
is \(\mathcal{N}_i = \{ j : (i,j)\in\mathcal{E}
\}\), and the associated Laplacian matrix
\(\mathbf{L}\in\mathbb{R}^{N\times N}\) is
defined by \([\mathbf{L}]_{ii} =
|\mathcal{N}_i|\) and \([\mathbf{L}]_{ij} = -1\)
if \((i,j)\in\mathcal{E}\), and \(0\) otherwise.
The sensing matrix \(\mathbf{A}\) and the
measurement vector \(\mathbf{b}\) are partitioned
row-wise across the agents, so that each agent
\(i\in\mathcal{V}\) holds a local data pair
\((\mathbf{A}_i,\mathbf{b}_i)\) with
\(\mathbf{A}_i\in\mathbb{R}^{m_i\times n}\),
\(\mathbf{b}_i\in\mathbb{R}^{m_i}\) and
\(\sum_{i=1}^N m_i = m\).  To enable fully
distributed computation without a central
coordinator, each agent maintains a local copy
\(\mathbf{x}_i\in\mathbb{R}^n\) of the global
decision variable \(\mathbf{x}\), and the
consensus constraint \(\mathbf{x}_1 =
\mathbf{x}_2 = \cdots = \mathbf{x}_N\) is
imposed. By introducing the above multi-agent
network, problem (\ref{eq:mollified_central}) can
be reformulated as the following distributed
optimization problem:
\begin{equation}
\label{eq:dist_mollified}
\begin{aligned}
\min_{\overline{\mathbf{x}}\in\mathbb{R}^{Nn}} \quad & \sum_{i=1}^{N} f_i(\mathbf{x}_i) \\
\text{s.t.} \quad & \mathbf{x}_1 = \cdots = \mathbf{x}_N,
\end{aligned}
\end{equation}
where \(f_i(\mathbf{x}_i) =
\frac{1}{2}\|\mathbf{A}_i\mathbf{x}_i -
\mathbf{b}_i\|^2 +
\frac{1}{N}P_{\lambda}^{\varepsilon}(\mathbf{x}_i)\),
and \(\overline{\mathbf{x}} =
\operatorname{col}(\mathbf{x}_1,\dots,\mathbf{x}_N)\).
Let \(\mathcal{C} = \{\overline{\mathbf{x}} \in
\mathbb{R}^{Nn} :
(\mathbf{L}\otimes\mathbf{I}_n)\overline{\mathbf{x}}
= \mathbf{0}\}\) denote the consensus subspace.
Observe that
\(\overline{\mathbf{x}}^{\mathsf{T}}(\mathbf{L}\otimes\mathbf{I}_n)\overline{\mathbf{x}}
= \frac{1}{2}\sum_{i=1}^N\sum_{j\in\mathcal{N}_i}
\|\mathbf{x}_i - \mathbf{x}_j\|^2 = 0\) for any
\(\overline{\mathbf{x}} \in \mathcal{C}\). Hence,
problem (\ref{eq:dist_mollified}) can be
equivalently reformulated as:
\begin{equation}
\label{eq:reformulated}
\begin{aligned}
\min_{\overline{\mathbf{x}} \in \mathbb{R}^{Nn}} \quad &  \sum_{i=1}^{N} f_i(\mathbf{x}_i) + \frac{1}{2}\overline{\mathbf{x}}^{\mathsf{T}}(\mathbf{L}\otimes\mathbf{I}_n)\overline{\mathbf{x}} \\
\text{s.t.} \quad & (\mathbf{L}\otimes\mathbf{I}_n)\overline{\mathbf{x}} = \mathbf{0}.
\end{aligned}
\end{equation}
The added quadratic term induces local convexity
along directions orthogonal to \(\mathcal{C}\),
which stabilizes the continuous-time dynamics to
be designed later.

By introducing Lagrangian multipliers to relax
the constraints~\cite{antil2026duality}, the
first-order necessary optimality conditions of
Problem (\ref{eq:reformulated}) can be formulated
as the following saddle-point problem:
\begin{equation}
\label{eq:saddle_point}
\min_{\overline{\mathbf{x}} \in \mathbb{R}^{Nn}} \; \max_{\overline{\bm{\omega}} \in \mathbb{R}^{Nn}} \; \mathcal{L}_{\varepsilon}(\overline{\mathbf{x}}, \overline{\bm{\omega}}),
\end{equation}
where \(\overline{\bm{\omega}} =
\operatorname{col}(\bm{\omega}_1,\dots,\bm{\omega}_N)\in\mathbb{R}^{Nn}\)
are the Lagrangian multipliers and
$\mathcal{L}_{\varepsilon}(\overline{\mathbf{x}},
\overline{\bm{\omega}}) = \sum_{i=1}^{N} \bigl(
f_i(\mathbf{x}_i) \bigr) +
\frac{1}{2}\overline{\mathbf{x}}^{\mathsf{T}}(\mathbf{L}\otimes\mathbf{I}_n)\overline{\mathbf{x}}
+
\overline{\bm{\omega}}^{\mathsf{T}}(\mathbf{L}\otimes\mathbf{I}_n)\overline{\mathbf{x}}$.

To solve \eqref{eq:saddle_point}, we develop a
distributed collaborative neurodynamic multi-agent
system that incorporates the mollified SCAD
penalty and exploits two-timescale singular
perturbation structure.  To each agent \(i\), we
assign a local momentum variable
\(\mathbf{m}_i\in\mathbb{R}^n\) and propose the
continuous-time evolution of the triplet
\((\mathbf{m}_i, \mathbf{x}_i, \bm{\omega}_i)\) by
\begin{equation}
\label{eq:dynamics_local}
\left\{
\begin{aligned}
\epsilon_1 \dot{\mathbf{m}}_i &= -\eta \mathbf{m}_i - \mathbf{A}_i^{\mathsf{T}}(\mathbf{A}_i\mathbf{x}_i - \mathbf{b}_i) - \frac{1}{N}\nabla P_{\lambda}^{\varepsilon}(\mathbf{x}_i) \\
&\quad - \sum_{j \in \mathcal{N}_i} \big( \bm{\omega}_i - \bm{\omega}_j + \mathbf{x}_i - \mathbf{x}_j \big), \\[4pt]
\epsilon_2 \dot{\mathbf{x}}_i &= \mathbf{m}_i, \\[4pt]
\epsilon_1 \dot{\bm{\omega}}_i &= \sum_{j \in \mathcal{N}_i} \big( \mathbf{m}_i - \mathbf{m}_j + \mathbf{x}_i - \mathbf{x}_j \big),
\end{aligned}
\right.
\end{equation}
where \(\epsilon_1,\epsilon_2\) are time-scale
parameters with $0<\epsilon_1\ll\epsilon_2$ and
\(\eta\) is the positive damping coefficient.
Multi-agent system~\eqref{eq:dynamics_local}
represents a three-layer recurrent neurodynamic network
at each agent. The momentum layer is a leaky
integrator that receives the local gradient and
neighbor coupling, and feeds into the state
layer. The state layer is a slow pure integrator
that accumulates the momentum signal, smoothing
out high-frequency fluctuations. The dual layer
is a fast, leakage-free integrator
that accumulates consensus errors and feeds back
to the momentum layer, closing the recurrent
loop.

Let \(\overline{\mathbf{m}} =
\operatorname{col}(\mathbf{m}_1,\dots,\mathbf{m}_N)\),
\(\overline{\bm{\omega}} =
\operatorname{col}(\bm{\omega}_1,\dots,\bm{\omega}_N)\),
\(\overline{\mathbf{A}} =
\operatorname{diag}(\mathbf{A}_1,\dots,\mathbf{A}_N)\),
\(\overline{\mathbf{b}} =
\operatorname{col}(\mathbf{b}_1,\dots,\mathbf{b}_N)\),
\(
\nabla \overline{P}_{\lambda}^{\varepsilon}(\overline{\mathbf{x}}) := \operatorname{col}\!\big(\tfrac{1}{N}\nabla P_{\lambda}^{\varepsilon}(\mathbf{x}_1),\dots,\tfrac{1}{N}\nabla P_{\lambda}^{\varepsilon}(\mathbf{x}_N)\big).
\)
Using the Laplacian identity
\(\sum_{j\in\mathcal{N}_i}(\mathbf{z}_i-\mathbf{z}_j)
=
[(\mathbf{L}\otimes\mathbf{I}_n)\overline{\mathbf{z}}]_i\),
the system~\eqref{eq:dynamics_local} can be
written in the compact form
\begin{equation}
\label{eq:dynamics_global}
\left\{
\begin{aligned}
\epsilon_1 \dot{\overline{\mathbf{m}}} &= -\eta \overline{\mathbf{m}} - \overline{\mathbf{A}}^{\mathsf{T}}(\overline{\mathbf{A}}\overline{\mathbf{x}} - \overline{\mathbf{b}}) \\
&\quad - \nabla \overline{P}_{\lambda}^{\varepsilon}(\overline{\mathbf{x}}) - (\mathbf{L}\otimes\mathbf{I}_n)(\overline{\bm{\omega}} + \overline{\mathbf{x}}), \\[4pt]
\epsilon_2 \dot{\overline{\mathbf{x}}} &= \overline{\mathbf{m}}, \\[4pt]
\epsilon_1 \dot{\overline{\bm{\omega}}} &= (\mathbf{L}\otimes\mathbf{I}_n)(\overline{\mathbf{m}} + \overline{\mathbf{x}}).
\end{aligned}
\right.
\end{equation}

We now characterize the equilibria of the collaborative neurodynamic multi-agent system~\eqref{eq:dynamics_global} and their relationship to the stationary points of the mollified problem.

\begin{theorem}
\label{thm:equilibrium} Let
\(\overline{\mathbf{z}} =
\operatorname{col}(\overline{\mathbf{m}},
\overline{\mathbf{x}}, \overline{\bm{\omega}})
\in \mathbb{R}^{3Nn}\) denote the state vector.
Assume that the communication graph
\(\mathcal{G}\) is connected \textbf{(C1)}.  Then
\(\overline{\mathbf{z}}^* =
\operatorname{col}(\overline{\mathbf{m}}^*,
\overline{\mathbf{x}}^*,
\overline{\bm{\omega}}^*)\) is an equilibrium of
system~\eqref{eq:dynamics_global} if and only if
the following conditions hold
\begin{enumerate}
    \item \(\overline{\mathbf{m}}^* =
        \mathbf{0}\).
    \item \(\overline{\mathbf{x}}^* =
        \mathbf{1}_N \otimes \mathbf{x}^*\) for
        some \(\mathbf{x}^* \in \mathbb{R}^n\),
        and this \(\mathbf{x}^*\) satisfies
        \(\sum_{i=1}^N
        \mathbf{A}_i^{\mathsf{T}}(\mathbf{A}_i\mathbf{x}^*
        - \mathbf{b}_i) + \nabla
        P_{\lambda}^{\varepsilon}(\mathbf{x}^*)
        = \mathbf{0}\).
    \item \(\overline{\bm{\omega}}^*\)
        satisfies
        \((\mathbf{L}\otimes\mathbf{I}_n)\overline{\bm{\omega}}^*
        = -
        \overline{\mathbf{A}}^{\mathsf{T}}(\overline{\mathbf{A}}\overline{\mathbf{x}}^*
        - \overline{\mathbf{b}}) - \nabla
        \overline{P}_{\lambda}^{\varepsilon}(\overline{\mathbf{x}}^*)\).
\end{enumerate}
\end{theorem}

\begin{proof}

\textbf{Necessity (\(\Rightarrow\)).}  Assume
\(\overline{\mathbf{z}}^*\) is an equilibrium.
Substituting \(\overline{\mathbf{z}}^*\) into the
second and third equations
of~\eqref{eq:dynamics_global} gives
\(\overline{\mathbf{m}}^* = \mathbf{0}\)
 and
\((\mathbf{L}\otimes\mathbf{I}_n)\overline{\mathbf{x}}^*
= \mathbf{0}\).  By connectedness of
\(\mathcal{G}\),
\(\operatorname{null}(\mathbf{L}) =
\operatorname{span}\{\mathbf{1}_N\}\); hence
there exists \(\mathbf{x}^* \in \mathbb{R}^n\)
such that \(\overline{\mathbf{x}}^* =
\mathbf{1}_N \otimes \mathbf{x}^*\).  With
\(\overline{\mathbf{m}}^* = \mathbf{0}\) and
\((\mathbf{L}\otimes\mathbf{I}_n)\overline{\mathbf{x}}^*
= \mathbf{0}\), the first equation
of~\eqref{eq:dynamics_global} at equilibrium
reduces to
\begin{equation} \label{eq:kkt_global}
\overline{\mathbf{A}}^{\mathsf{T}}(\overline{\mathbf{A}}\overline{\mathbf{x}}^* - \overline{\mathbf{b}}) + \nabla \overline{P}_{\lambda}^{\varepsilon}(\overline{\mathbf{x}}^*) + (\mathbf{L}\otimes\mathbf{I}_n)\overline{\bm{\omega}}^* = \mathbf{0}.
\end{equation}
Left-multiplying~\eqref{eq:kkt_global} by
\(\mathbf{1}_N^{\mathsf{T}} \otimes
\mathbf{I}_n\) and using the fact that
\(\mathbf{1}_N^{\mathsf{T}}\mathbf{L} =
\mathbf{0}^{\mathsf{T}}\), the block-diagonal
structure of \(\overline{\mathbf{A}}\), and the
identity \( (\mathbf{1}_N^{\mathsf{T}} \otimes
\mathbf{I}_n) \nabla
\overline{P}_{\lambda}^{\varepsilon}(\overline{\mathbf{x}})
= \nabla P_{\lambda}^{\varepsilon}(\mathbf{x}) \)
whenever \(\overline{\mathbf{x}} = \mathbf{1}_N
\otimes \mathbf{x}\), we obtain precisely
item~(2). Item~(3) follows directly
from~\eqref{eq:kkt_global}.

\textbf{Sufficiency (\(\Leftarrow\)).}
Suppose \(\overline{\mathbf{z}}^*\) satisfies items~(1)--(3).
Item~(1) immediately fulfills the second equation of~\eqref{eq:dynamics_global}.
Item~(2) gives \(\overline{\mathbf{x}}^* = \mathbf{1}_N \otimes \mathbf{x}^*\), which by connectivity implies \((\mathbf{L}\otimes\mathbf{I}_n)\overline{\mathbf{x}}^* = \mathbf{0}\); together with \(\overline{\mathbf{m}}^* = \mathbf{0}\) this satisfies the third equation.
Substituting \(\overline{\mathbf{m}}^* = \mathbf{0}\) and \((\mathbf{L}\otimes\mathbf{I}_n)\overline{\mathbf{x}}^* = \mathbf{0}\) into the first equation yields the left-hand side \(\overline{\mathbf{A}}^{\mathsf{T}}(\overline{\mathbf{A}}\overline{\mathbf{x}}^* - \overline{\mathbf{b}}) + \nabla \overline{P}_{\lambda}^{\varepsilon}(\overline{\mathbf{x}}^*) + (\mathbf{L}\otimes\mathbf{I}_n)\overline{\bm{\omega}}^*\), which vanishes by item~(3).
Thus all three equations are satisfied at \(\overline{\mathbf{z}}^*\).
\end{proof}

\begin{remark}
\label{rem:equilibrium_stationarity}
Theorem~\ref{thm:equilibrium} characterizes the
equilibria of~\eqref{eq:dynamics_global} via
three equivalent conditions. For any equilibrium
\(\overline{\mathbf{z}}^* =
\operatorname{col}(\overline{\mathbf{m}}^*,
\overline{\mathbf{x}}^*,
\overline{\bm{\omega}}^*)\), items~(2) and~(3)
imply that \((\overline{\mathbf{x}}^*,
\overline{\bm{\omega}}^*)\) is a KKT point of the
saddle-point problem~\eqref{eq:saddle_point}. By
consensus (item~(2)),
\(\overline{\mathbf{x}}^* = \mathbf{1}_N \otimes
\mathbf{x}^*\) for some \(\mathbf{x}^* \in
\mathbb{R}^n\); substituting this structure into
the primal KKT condition recovers exactly the
stationarity condition of the central smoothed
problem~\eqref{eq:mollified_central}. Thus
\(\mathbf{x}^*\) is a stationary point of the
smooth surrogate objective, establishing a
one-to-one correspondence between the equilibria
of the collaborative neurodynamic multi-agent system and the
stationary points of the smoothed problem of
interest.
\end{remark}

Having characterized the equilibria, we next
establish local optimality conditions under which
an equilibrium of the distributed system
corresponds to a strict local minimum of the
centralized problem. Let \(I_1\) and \(I_0\)
denote the index sets of nonzero and zero entries
of \(\mathbf{x}^*\), respectively. Let
\begin{equation*}
\mathbf{H}_{I_{1,1}} := \Bigl[\frac{1}{N}\sum_{i=1}^N \mathbf{A}_i^{\mathsf{T}}\mathbf{A}_i\Bigr]_{I_1,I_1}
\end{equation*}
denote the principal submatrix of the averaged
Gram matrix indexed by \(I_1\), and let
\(\mathbf{H}_{I_{1,0}}\) denote the corresponding
cross block between \(I_1\) and \(I_0\).

\begin{theorem}
\label{thm:local_optimality} Suppose there exists
\(\delta>0\) such that \(|x_j^*| \ge a\lambda +
\delta\) for all \(j \in I_1\) and \(x_j^* = 0\)
for all \(j \in I_0\) \textbf{(C2)}. Let
\(\lambda_{I_1}:=\lambda_{\min}(\mathbf{H}_{I_1})>0\).
If the mollification parameter satisfies
\begin{equation}
\varepsilon \le \min\left\{ \delta,\; \frac{2\lambda c \lambda_{I_1}}{e\left( 4N\lambda_{I_1}\|\mathbf{H}_{I_{1,0}}\|^2 + \frac{\lambda_{I_1}}{a-1} \right)} \right\},
\label{eq:eps_condition}
\end{equation}
\textbf{(C3)}, then \(\mathbf{x}^*\) is a strict
local minimum of the central mollified
problem~\eqref{eq:mollified_central}, and
\(\overline{\mathbf{x}}^* = \mathbf{1}_N \otimes
\mathbf{x}^*\) is a strict local minimum of the
distributed optimization
problem~\eqref{eq:dist_mollified}.
\end{theorem}

\begin{proof}
Let \(F(\mathbf{x}) = \frac{1}{N}\sum_{i=1}^N
f_i(\mathbf{x})\). Its Hessian at
\(\mathbf{x}^*\) is
\begin{equation*}
\nabla^2 F(\mathbf{x}^*) = \frac{1}{N}\sum_{i=1}^N \mathbf{A}_i^{\mathsf{T}}\mathbf{A}_i + \frac{1}{N}\operatorname{diag}\bigl(q_{\varepsilon}'(x_1^*),\dots,q_{\varepsilon}'(x_n^*)\bigr).
\end{equation*}
For \(j\in I_1\), \(\varepsilon<\delta\) gives
\(q_{\varepsilon}'(x_j^*)=0\). For \(j\in I_0\),
\(x_j^*=0\). Let \( p_{\lambda,\text{reg}}'' \)
denote the pointwise second derivative of the
SCAD penalty function outside the set of points
\( \{0, \pm\lambda, \pm a\lambda\} \). A
straightforward convolution smoothing computation
yields
\begin{equation}
q_\varepsilon'(0) =
(\rho_\varepsilon * p_{\lambda,\text{reg}}'')(0) +
2\lambda\rho_\varepsilon(0) \ge -\frac{1}{a-1} +
\frac{2c\lambda}{e\varepsilon}.
\label{eq:q_derivative_lower_bound}
\end{equation}
Let \(\mathbf{D} = d_0\mathbf{I}_{I_0}\) and
\(d_0  = q_{\varepsilon}'(0)/N \). Combining~%
(\ref{eq:eps_condition}) and~%
(\ref{eq:q_derivative_lower_bound}) yields \( d_0
>
4\|\mathbf{H}_{I_{1,0}}\|^2/\lambda_{I_1}
> 0\).
Partitioning the averaged Gram matrix according
to \(I_1\) and \(I_0\) yields
\begin{equation*}
\nabla^2 F(\mathbf{x}^*) =
\begin{pmatrix}
\mathbf{H}_{I_1} & \mathbf{H}_{I_{1,0}}\\[4pt]
\mathbf{H}_{I_{1,0}}^{\mathsf{T}} & \mathbf{H}_{I_0} + \mathbf{D}
\end{pmatrix}.
\end{equation*}
For any
\(\mathbf{v}=(\mathbf{v}_{I_1},\mathbf{v}_{I_0})\),
applying Schur complement and Young's inequality
yields
\begin{align*}
\mathbf{v}^{\mathsf{T}} \nabla^2 F(\mathbf{x}^*) \mathbf{v}
&\ge \lambda_{I_1} \|\mathbf{v}_{I_1}\|^2 + d_0 \|\mathbf{v}_{I_0}\|^2 \\
&\quad - 2\|\mathbf{H}_{I_{1,0}}\| \|\mathbf{v}_{I_1}\| \|\mathbf{v}_{I_0}\| \\
&\ge \frac{\lambda_{I_1}}{2}\|\mathbf{v}_{I_1}\|^2
   + \bigl(d_0 - 2\|\mathbf{H}_{I_{1,0}}\|^2/\lambda_{I_1}\bigr)\|\mathbf{v}_{I_0}\|^2 \\
&\ge \mu_0 \|\mathbf{v}\|^2,
\end{align*}
where \(\mu_0
=\min\bigl(\lambda_{I_1}/2,\,2\|\mathbf{H}_{I_{1,0}}\|^2/\lambda_{I_1}\bigr)>0\)
. Hence, \(\nabla^2 F(\mathbf{x}^*)\succeq \mu_0
\mathbf{I}\).

Since \(F(\mathbf{x})\) is \(C^2\), there exists
a convex neighborhood \(\mathcal{X}\) of
\(\mathbf{x}^*\) such that \(\nabla^2
F(\mathbf{x}) \succeq \frac{\mu_0}{2}\mathbf{I}\)
for all \(\mathbf{x}\in\mathcal{X}\); hence
\(F(\mathbf{x})\) is strongly convex on
\(\mathcal{X}\) and \(\mathbf{x}^*\) is its
unique minimizer on \(\mathcal{X}\). In
particular, \(\mathbf{x}^*\) is a strict local
minimum of \(F(\mathbf{x})\), and because the
central objective in~\eqref{eq:mollified_central}
is a positive multiple of \(F(\mathbf{x})\), it
is also a strict local minimum of that problem.
For the distributed
problem~\eqref{eq:dist_mollified}, feasibility
forces \(\overline{\mathbf{x}} = \mathbf{1}_N
\otimes \mathbf{x}\) and the objective equals \(N
F(\mathbf{x})\).  The norm relation
\(\|\overline{\mathbf{x}} -
\overline{\mathbf{x}}^*\| = \sqrt{N}\|\mathbf{x}
- \mathbf{x}^*\|\) ensures that points
sufficiently close to \(\overline{\mathbf{x}}^*\)
project to \(\mathbf{x} \in \mathcal{X}\); on
this set, \(N F(\mathbf{x}) \ge N
F(\mathbf{x}^*)\) with equality only at
\(\mathbf{x} = \mathbf{x}^*\).  Thus
\(\overline{\mathbf{x}}^*\) is a strict local
minimum of the distributed problem under the
consensus constraint.
\end{proof}

The local optimality established in
Theorem~\ref{thm:local_optimality} implies the
following gradient growth condition around
\(\mathbf{x}^*\).

\begin{corollary}
\label{cor:gradient_growth} Under conditions
\textbf{(C2)}--\textbf{(C3)}, there exists a
neighborhood \(\mathcal{X}\) of \(\mathbf{x}^*\)
and a constant
\[
\mu = \frac14 \min\!\left(\lambda_{I_1},\; \frac{4\|\mathbf{H}_{I_{1,0}}\|^2}{\lambda_{I_1}}\right) > 0
\]
such that
\[
\|\nabla F(\mathbf{x})\| \ge \mu \|\mathbf{x} - \mathbf{x}^*\|, \quad \forall \mathbf{x}\in\mathcal{X}.
\]
\end{corollary}

\begin{proof}
According to the proof of
Theorem~\ref{thm:local_optimality},
\(F(\mathbf{x})\) is strongly convex on
\(\mathcal{X}\) with modulus \(\mu_0/2\). Since
\(\nabla F(\mathbf{x}^*) = \mathbf{0}\), the
standard strong convexity inequality gives
\[
\langle \nabla F(\mathbf{x}), \mathbf{x} - \mathbf{x}^* \rangle
\ge \frac{\mu_0}{2} \|\mathbf{x} - \mathbf{x}^*\|^2, \quad \mathbf{x}\in\mathcal{X}.
\]
Cauchy--Schwarz yields \(\|\nabla F(\mathbf{x})\|
\ge \frac{\mu_0}{2} \|\mathbf{x} -
\mathbf{x}^*\|\). Hence the claimed constant is
\(\mu = \mu_0/2\), which completes the proof.
\end{proof}

\begin{remark}
\label{rem:equilibrium_original_relation}
Theorem~\ref{thm:local_optimality} shows that
\(\overline{\mathbf{x}}^* = \mathbf{1}_N \otimes
\mathbf{x}^*\) is a strict local minimum of the
smoothed distributed
problem~\eqref{eq:dist_mollified}. Since the
mollified penalty deviates from the original SCAD
subgradient only on a set of measure
\(O(\varepsilon)\), this point automatically
satisfies the Clarke stationarity condition
for~\eqref{eq:sparse_lsq}. As
\(\varepsilon\to0^+\) the discrepancy set
vanishes; consequently, the equilibrium converges
to a local solution of the original non-smooth
SCAD problem.
\end{remark}

\begin{remark}
\label{rem:unbiasedness_preservation}
Under conditions \textbf{(C2)}--\textbf{(C3)},
the structure-preserving mollification
preserves the unbiasedness of the SCAD
penalty on large coefficients. If
\(|x_j^*| \ge a\lambda + \delta\) for
\(j\in I_1\) and \(\varepsilon < \delta\),
then the mollifier support
\([-\varepsilon,\varepsilon]\) does not
reach the large-coefficient region, so the
mollified derivative satisfies
\(q_\varepsilon(x_j^*) = p_\lambda'(x_j^*) = 0\)
and the mollified penalty coincides
pointwise with the original SCAD penalty
on a neighborhood of \(\mathbf{x}^*_{I_1}\).
The asymptotic unbiasedness and normality
of Fan and Li~\cite{fan2001variable}, which depend only on the
penalty behavior in this neighborhood, are
therefore preserved. Note, however, that
the sparsity selection component of the
oracle property is not preserved by the
mollification alone; see
Remark~\ref{rem:sparsity_pso}.
\end{remark}

\section{Convergence Analysis}
\label{sec:convergence}

We now establish convergence guarantees for the collaborative neurodynamic multi-agent system. The analysis proceeds via singular perturbation theory~\cite{khalil2002nonlinear, kokotovic1986singular}. We first construct a reduced fast-slow decomposition, then derive trajectory bounds through Lyapunov synthesis.

Define the mapping \(\mathbf{F}_{\rm map} : \mathbb{R}^{3Nn} \to \mathbb{R}^{2Nn}\) by
\[
\mathbf{F}_{\rm map}(\overline{\mathbf{z}}) =
\begin{bmatrix}
-\eta \overline{\mathbf{m}} - (\mathbf{L}\otimes\mathbf{I}_n)\overline{\bm{\omega}} - \bm{\Phi}(\overline{\mathbf{x}}) \\[3pt]
(\mathbf{L}\otimes\mathbf{I}_n)(\overline{\mathbf{m}} + \overline{\mathbf{x}})
\end{bmatrix},
\]
where
$
\bm{\Phi}(\overline{\mathbf{x}}) = \overline{\mathbf{A}}^{\mathsf{T}}(\overline{\mathbf{A}}\overline{\mathbf{x}} - \overline{\mathbf{b}}) + \nabla\overline{P}_{\lambda}^{\varepsilon}(\overline{\mathbf{x}}) + (\mathbf{L}\otimes\mathbf{I}_n)\overline{\mathbf{x}}.
$
Let \(\overline{\mathbf{p}} =
\operatorname{col}(\overline{\mathbf{m}},
\overline{\bm{\omega}})\). We propose to express
\(\mathbf{F}_{\rm map}(\overline{\mathbf{z}})\) via
the affine decomposition
\[
\mathbf{F}_{\rm map}(\overline{\mathbf{z}}) = \mathbf{G}_{\rm c}\overline{\mathbf{p}} + \mathbf{G}(\overline{\mathbf{x}}),
\]
with
\begin{equation*}
\mathbf{G}_{\rm c} = \begin{bmatrix}
-\eta \mathbf{I}_{Nn} & -\mathbf{L}\otimes\mathbf{I}_n \\[3pt]
\mathbf{L}\otimes\mathbf{I}_n & \mathbf{0}_{Nn}
\end{bmatrix}, \quad
\mathbf{G}(\overline{\mathbf{x}}) = \begin{bmatrix}
-\bm{\Phi}(\overline{\mathbf{x}}) \\[3pt]
(\mathbf{L}\otimes\mathbf{I}_n)\overline{\mathbf{x}}
\end{bmatrix}.
\end{equation*}
Let
\begin{equation*}
  \mathcal{B}_{\overline{\mathbf{x}}^*} = \bigl\{ \overline{\mathbf{x}} \in \mathbb{R}^{Nn} : \|\overline{\mathbf{x}} - \overline{\mathbf{x}}^*\| \le r \bigr\}.
\end{equation*}
where \(r > 0\) is chosen sufficiently small so
that \(\mathcal{B}_{\overline{\mathbf{x}}^*}
\subset \{ \mathbf{1}_N \otimes \mathbf{x} :
\mathbf{x} \in \mathcal{X} \}\); such a radius
exists because \(\mathcal{X}\) is a neighborhood
of \(\mathbf{x}^*\). Let \(l_1 = \max_i
\|\mathbf{A}_i^{\mathsf{T}}\mathbf{A}_i\|\).

We first characterize the fast-variable
dynamics during the boundary layer, then
establish local ultimate boundedness of the
full system trajectories around the
equilibrium.

\begin{proposition}
\label{prop:boundary_layer}
Under conditions \textbf{(C1)}--\textbf{(C4)},
the fast variable \(\overline{\mathbf{p}} =
\operatorname{col}(\overline{\mathbf{m}},
\overline{\boldsymbol{\omega}})\) converges
exponentially to the quasi-steady-state
manifold during the boundary layer
\([0,T]\). On the fast time scale
\(\tau = t/\epsilon_1\), the tracking error
\(\mathbf{e}(\tau) = \overline{\mathbf{p}}(\tau) -
\overline{\mathbf{p}}_{\rm qss}\) satisfies
\[
\|\mathbf{e}(\tau)\| \le
K\,\|\mathbf{e}(0)\|\,e^{-\alpha \tau},
\]
where \(\overline{\mathbf{p}}_{\rm qss} =
-\mathbf{G}_{\rm c}^{-1}\mathbf{r}(\overline{\mathbf{x}})\),
\(\alpha = \min_k |\operatorname{Re}(\lambda_k)| > 0\)
for \(\{\lambda_k\}\) the spectrum of the
Hurwitz matrix \(\mathbf{G}_{\rm c}\), and
\(K = \kappa(\mathbf{V})\) is the condition
number of the eigenvector matrix
\(\mathbf{V}\). The boundary layer duration
\(T = 5\epsilon_1/\alpha\) ensures
\(\|\mathbf{e}(T/\epsilon_1)\| \le
0.7\%\,K\,\|\mathbf{e}(0)\|\).
\end{proposition}

\begin{proof}
On the fast time scale \(\tau = t/\epsilon_1\),
the slow variable evolves at rate
\(d\overline{\mathbf{x}}/d\tau =
\epsilon_1\dot{\overline{\mathbf{x}}} =
O(\epsilon_1)\); by Tikhonov's
theorem~\cite{khalil2002nonlinear}
it is treated as frozen during the boundary
layer, so the fast subsystem reduces to the
affine ODE
\[
\frac{d\overline{\mathbf{p}}}{d\tau} =
\mathbf{F}_{\rm map}(\overline{\mathbf{z}}) =
\mathbf{G}_{\rm c}\,\overline{\mathbf{p}} +
\mathbf{r}(\overline{\mathbf{x}}),
\]
with unique quasi-steady-state
\(\overline{\mathbf{p}}_{\rm qss} =
-\mathbf{G}_{\rm c}^{-1}\mathbf{r}(\overline{\mathbf{x}})\).
Because \(\mathbf{F}_{\rm map}\) is affine in
\(\overline{\mathbf{p}}\), the linearization
about \(\overline{\mathbf{p}}_{\rm qss}\) is
exact rather than first-order: setting
\(\mathbf{e} = \overline{\mathbf{p}} -
\overline{\mathbf{p}}_{\rm qss}\) gives the
homogeneous linear system
\[
\frac{d\mathbf{e}}{d\tau} = \mathbf{G}_{\rm c}\,\mathbf{e},
\qquad
\mathbf{e}(0) = \overline{\mathbf{p}}(0) -
\overline{\mathbf{p}}_{\rm qss},
\]
with solution
\(\mathbf{e}(\tau) = e^{\mathbf{G}_{\rm c}\tau}\,\mathbf{e}(0)\).
Let \(\mathbf{G}_{\rm c} = \mathbf{V}\mathbf{J}\mathbf{V}^{-1}\)
be the Jordan decomposition, with
\(\operatorname{Re}(\lambda_k) \le -\alpha\)
for every eigenvalue. Then
\[
\|e^{\mathbf{G}_{\rm c}\tau}\|
\le \|\mathbf{V}\|\,\|\mathbf{V}^{-1}\|\,\|e^{\mathbf{J}\tau}\|
\le \kappa(\mathbf{V})\,e^{-\alpha \tau},
\]
where the last step applies the standard
Jordan-block exponential bound for Hurwitz
matrices. Setting \(K = \kappa(\mathbf{V})\)
yields the claimed estimate. The choice
\(T = 5/\alpha\) in fast time, i.e.,
\(T = 5\epsilon_1/\alpha\) in slow time,
gives \(\|\mathbf{e}(T/\epsilon_1)\| \le
K\|\mathbf{e}(0)\|e^{-5} \le 0.7\%
K\|\mathbf{e}(0)\|\).
\end{proof}

\begin{lemma}
\label{lem:main} Suppose conditions
\textbf{(C1)}--\textbf{(C3)} hold and
$\overline{\mathbf{x}}(0) \in
\mathcal{B}_{\overline{\mathbf{x}}^*}$ \textbf{(C4)}. Then,
when $\epsilon_1 = 0$, the exact
quasi-steady-state inequality
\[
\|\overline{\mathbf{x}} - \overline{\mathbf{x}}^*\| \le c_0 \|\overline{\mathbf{m}}\|
\]
holds. Moreover, there exist constants
\(\epsilon_1^*>0\), \(c_1>0\) and  a finite time
\(T \ge 0\) such that the inequality
\[
\|\overline{\mathbf{x}}(t) - \overline{\mathbf{x}}^*\| \le c_0 \|\overline{\mathbf{m}}(t)\| + c_1 \epsilon_1
\]
holds for every \(\epsilon_1 \in
(0,\epsilon_1^*)\) and \(t \ge T\), where \(c_0 =
\sqrt{1 + (\eta + l_1 + l_\varepsilon/N)^2 /
\mu^2}\).
\end{lemma}
\begin{proof}

According to singular perturbation theory
\cite{khalil2002nonlinear, wang2022two}, the fast
subsystem in (\ref{eq:dynamics_global}) can be
expressed as
\[
\epsilon_1 \dot{\overline{\mathbf{p}}} = \mathbf{F}_{\rm map}(\overline{\mathbf{z}}).
\]
Since
\(\mathbf{G}_{\rm c}\) is Hurwitz on the subspace
orthogonal to consensus, it satisfies the
exponential stability condition of the Tikhonov
theorem. Consequently, letting \(\epsilon_1 = 0\)
yields the algebraic equations
$\mathbf{F}_{\rm map}(\overline{\mathbf{z}}) = 0$ on
the  quasi-steady-state manifold. Let
\(\bar{\mathbf{m}} = \frac{1}{N}\sum_{i=1}^N
\mathbf{m}_i\) and \(\bar{\mathbf{x}} =
\frac{1}{N}\sum_{i=1}^N \mathbf{x}_i\). It
follows from a straightforward argument that
\begin{align}
  \mathbf{x}_i - \bar{\mathbf{x}} &= -(\mathbf{m}_i - \bar{\mathbf{m}}), \label{eq:manifold_id1} \\
  \frac{1}{N}\sum_{j=1}^N \nabla f_j(\mathbf{x}_j) &= -\eta\,\bar{\mathbf{m}}. \label{eq:manifold_id2}
\end{align}
holds on the manifold. It follows from
Lemma~\ref{lem:lip} that \(\nabla F(\bar{\mathbf{x}})
= \frac{1}{N}\sum_{i=1}^N \nabla
f_i(\bar{\mathbf{x}})\) and each \(\nabla
f_i(x_i)\) is Lipschitz continuous with constant
at most \(l_1 + l_\varepsilon/N\). From
(\ref{eq:manifold_id1}) and
(\ref{eq:manifold_id2}) together with the
Cauchy--Schwarz inequality, it follows that
\begin{equation}
  \begin{aligned}
    \|\nabla F(\bar{\mathbf{x}})\|
    &\le \Bigl\|\frac{1}{N}\sum_{i=1}^N \nabla f_i(\mathbf{x}_i)\Bigr\| \\
    &\quad + \frac{1}{N}\sum_{i=1}^N \|\nabla f_i(\mathbf{x}_i) - \nabla f_i(\bar{\mathbf{x}})\| \\
    &\le \eta\,\|\bar{\mathbf{m}}\| + \tfrac{l_1 + l_\varepsilon/N}{N}\sum_{i=1}^N \|\mathbf{x}_i - \bar{\mathbf{x}}\| \\
    &= \eta\,\|\bar{\mathbf{m}}\| + \tfrac{l_1 + l_\varepsilon/N}{N}\sum_{i=1}^N \|\mathbf{m}_i - \bar{\mathbf{m}}\|  \\
    & \le (\eta + l_1 + l_\varepsilon/N)\|\overline{\mathbf{m}}\|/\sqrt{N}.
  \end{aligned}
  \label{eq:grad_estimate}
\end{equation}
Under \textbf{(C1)}--\textbf{(C4)},
Corollary~\ref{cor:gradient_growth} and
(\ref{eq:grad_estimate}), we obtain
\begin{equation}
  \|\bar{\mathbf{x}} - \mathbf{x}^*\| \le \frac{\eta + l_1 + l_\varepsilon/N}{\mu\sqrt{N}}\,\|\overline{\mathbf{m}}\|.
  \label{eq:mean_error_bound}
\end{equation}
From ~\eqref{eq:mean_error_bound} and
~\eqref{eq:manifold_id1}, we have
\begin{equation*}
  \begin{aligned}
    \|\overline{\mathbf{x}} - \overline{\mathbf{x}}^*\|^2
    &= \sum_{i=1}^N \|\mathbf{x}_i - \bar{\mathbf{x}}\|^2 + N\|\bar{\mathbf{x}} - \mathbf{x}^*\|^2 \\
    &= \sum_{i=1}^N \|\mathbf{m}_i - \bar{\mathbf{m}}\|^2 + N\|\bar{\mathbf{x}} - \mathbf{x}^*\|^2 \\
    &\le \|\overline{\mathbf{m}}\|^2 + N \cdot \frac{(\eta + l_1 + l_\varepsilon/N)^2}{\mu^2 N}\,\|\overline{\mathbf{m}}\|^2 \\
    &= \Bigl(1 + \frac{(\eta + l_1 + l_\varepsilon/N)^2}{\mu^2}\Bigr)\,\|\overline{\mathbf{m}}\|^2.
  \end{aligned}
\end{equation*}
Taking square roots yields the exact manifold
bound \(\|\overline{\mathbf{x}} -
\overline{\mathbf{x}}^*\| \le
c_0\,\|\overline{\mathbf{m}}\|\).

By Tikhonov's theorem, for $\epsilon_1 > 0$,
there exist $\epsilon_1^* > 0$ and $c_1 > 0$
such that, for every $\epsilon_1 \in (0,
\epsilon_1^*)$, one can find a finite time
$T \ge 0$ for which the trajectory deviates
from the quasi-steady-state manifold by at most
$O(\epsilon_1)$, yielding the perturbed bound
\[
\|\overline{\mathbf{x}}(t) - \overline{\mathbf{x}}^*\| \le c_0\,\|\overline{\mathbf{m}}(t)\| + c_1\,\epsilon_1
\]
for all $t \ge T$.
\end{proof}

The next result sharpens the trajectory bound on the auxiliary variable \(\overline{\bm{\omega}}\) under the quasi-steady-state condition.

\begin{lemma}
\label{lem:omega_bound} Under conditions
\textbf{(C1)}--\textbf{(C4)}, the exact
quasi-steady-state inequality
\begin{equation*}
\|\overline{\boldsymbol{\omega}}(t) - \overline{\boldsymbol{\omega}}^*\| \le c_\omega \|\overline{\mathbf{m}}(t)\|
\end{equation*}
holds when \(\epsilon_1 = 0\), where the constant
\[
c_\omega = 1 + \frac{1}{\lambda_2(\mathbf{L})}\bigl(\eta + (l_1 + l_\varepsilon/N) c_0\bigr),
\]
\(\lambda_2(\mathbf{L})\) is the smallest nonzero
eigenvalue of the graph Laplacian \(\mathbf{L}\).
Moreover, there exist constants
\(\epsilon_1^*>0\), \(c_2>0, c_3>0\) and a finite
time \(T \ge 0\) such that the inequality
\begin{equation*}
\|\overline{\boldsymbol{\omega}}(t) - \overline{\boldsymbol{\omega}}^*\| \le c_\omega \|\overline{\mathbf{m}}(t)\| + (c_3 + c_\omega c_2)\epsilon_1
\end{equation*}
holds for every \(\epsilon_1 \in
(0,\epsilon_1^*)\) and \(t \ge T\).
\end{lemma}

\begin{proof}
Let \(\epsilon_1=0\)
in~\eqref{eq:dynamics_global} and arguing as in
Lemma~\ref{lem:main}, we obtain
$\mathbf{F}_{\rm map}(\overline{\mathbf{z}}) = 0$.
Further, there exists a vector \(\mathbf{v}\)
such that \(\overline{\mathbf{m}} +
\overline{\mathbf{x}} =
\mathbf{1}\otimes\mathbf{v}\).  Introducing the
zero-mean vectors \(\widetilde{\mathbf{m}} =
\overline{\mathbf{m}} -
\mathbf{1}\otimes\bar{\mathbf{m}},\
\widetilde{\mathbf{x}} = \overline{\mathbf{x}} -
\mathbf{1}\otimes\bar{\mathbf{x}}\), a simple
calculation yields \(\widetilde{\mathbf{x}} =
-\widetilde{\mathbf{m}}\).
 Observe that
\(\frac{1}{N}(\mathbf{1}^{\!\top}\otimes\mathbf{I}_n)\overline{\bm{\omega}}\)
is a conserved quantity by the third equation
in~\eqref{eq:dynamics_global}; consequently the
mean of \(\overline{\bm{\omega}}\) on the slow
manifold coincides with that of
\(\overline{\bm{\omega}}^*\), so
\(\overline{\bm{\omega}} -
\overline{\bm{\omega}}^*\) is zero-mean. The
equality $\mathbf{F}_{\rm map}(\overline{\mathbf{z}})
= 0$ gives
\[
(\mathbf{L}\otimes\mathbf{I}_n)\bigl(\overline{\bm{\omega}} - \overline{\bm{\omega}}^* + \overline{\mathbf{x}} - \overline{\mathbf{x}}^*\bigr)
= -\eta\overline{\mathbf{m}} - \bigl(\nabla F(\overline{\mathbf{x}}) - \nabla F(\overline{\mathbf{x}}^*)\bigr).
\]
Observe that \(\overline{\mathbf{x}} -
\overline{\mathbf{x}}^* = \widetilde{\mathbf{x}}
+ \mathbf{1}\otimes(\bar{\mathbf{x}} -
\mathbf{x}^*) = -\widetilde{\mathbf{m}} +
\mathbf{1}\otimes(\bar{\mathbf{x}} -
\mathbf{x}^*)\). Hence, we obtain
\[
(\mathbf{L}\otimes\mathbf{I}_n)\bigl(\overline{\bm{\omega}} - \overline{\bm{\omega}}^* - \widetilde{\mathbf{m}}\bigr)
= -\eta\overline{\mathbf{m}} - \bigl(\nabla F(\overline{\mathbf{x}}) - \nabla F(\overline{\mathbf{x}}^*)\bigr).
\]
The zero-mean property together with the graph's
connectivity yields
\[
\|\overline{\bm{\omega}} - \overline{\bm{\omega}}^* - \widetilde{\mathbf{m}}\|
\le \frac{1}{\lambda_2(\mathbf{L})}\Bigl(\eta\|\overline{\mathbf{m}}\| + \|\nabla F(\overline{\mathbf{x}}) - \nabla F(\overline{\mathbf{x}}^*)\|\Bigr).
\]
Using \(\|\overline{\bm{\omega}} -
\overline{\bm{\omega}}^*\| \le
\|\widetilde{\mathbf{m}}\| +
\|\overline{\bm{\omega}} -
\overline{\bm{\omega}}^* -
\widetilde{\mathbf{m}}\|\) and
\(\|\widetilde{\mathbf{m}}\| \le
\|\overline{\mathbf{m}}\|\), together with the
Lipschitz property of \(\nabla F(\mathbf{x})\)
 and the bound
\(\|\overline{\mathbf{x}} -
\overline{\mathbf{x}}^*\| \le
c_0\|\overline{\mathbf{m}}\|\) from
Lemma~\ref{lem:main}, we immediately obtain
\[
\|\overline{\bm{\omega}} - \overline{\bm{\omega}}^*\| \le
\Bigl(1 + \frac{\eta + (l_1 + l_\varepsilon/N) c_0}{\lambda_2(\mathbf{L})}\Bigr) \|\overline{\mathbf{m}}\|.
\]
This establishes the exact manifold inequality
with \(c_\omega\). For \(\epsilon_1>0\), similar
to Lemma~\ref{lem:main}, there exist constants
\(c_2,c_3>0\) such that, for \(t\ge T\),
$
\|\overline{\bm{\omega}}(t)-\overline{\bm{\omega}}^*\| \le c_\omega\|\overline{\mathbf{m}}(t)\| + (c_3+c_\omega c_2)\epsilon_1.
$
\end{proof}

Let $l_2 = \|\mathbf{G}_{\rm c}^{\mathsf T}\mathbf{G}_{\rm c}+\mathbf{I}\|$,
$l_3 = \sqrt{(l_1+l_\varepsilon/N)^2+\|\mathbf{L}\otimes\mathbf{I}_n\|^2}$,
and $\overline{\mathbf{p}}^* = \operatorname{col}(\overline{\mathbf{m}}^*,\overline{\bm{\omega}}^*)$.
Define
$
V_1(\overline{\mathbf{p}}) = \tfrac{1}{2}\|\mathbf{F}_{\rm map}(\overline{\mathbf{x}}^*,\overline{\mathbf{p}})\|^2 + \tfrac{1}{2}\|\overline{\mathbf{p}}-\overline{\mathbf{p}}^*\|^2.
$
The following corollary quantifies the residual cross-term arising from momentum-state interactions.
\begin{corollary}
\label{cor:local_bounds} Under conditions
\textbf{(C1)}--\textbf{(C4)},  the inequality
\[
\|\nabla V_1(\overline{\mathbf{p}})\|\,
\|\mathbf{F}_{\rm map}(\overline{\mathbf{z}})
-\mathbf{F}_{\rm map}(\overline{\mathbf{x}}^*,\overline{\mathbf{p}})\|
\le l_4\|\overline{\mathbf{m}}\|^2 + l_5\epsilon_1^2
\]  holds for every
$\epsilon_1\in(0,\epsilon_1^*)$ and $t\ge T$,
where $l_4 = (3/2)\,l_2 l_3 (1+c_\omega)\,c_0$
and $l_5 = l_2 l_3 \bigl\{ [(1+c_\omega)c_1 +
(c_3+c_\omega c_2)c_0]^2 / [2(1+c_\omega)c_0] +
(c_3+c_\omega c_2)c_1 \bigr\}$.
\end{corollary}

\begin{proof}
From Lemma~\ref{lem:omega_bound}, for every
\(\epsilon_1\in(0,\epsilon_1^*)\) and all \(t\ge
T\),
\[
\|\overline{\bm{\omega}}-\overline{\bm{\omega}}^*\|
\le c_\omega\|\overline{\mathbf{m}}\| + (c_3+c_\omega c_2)\epsilon_1.
\]
Hence the stacked variable
\(\overline{\mathbf{p}}=(\overline{\mathbf{m}},\overline{\bm{\omega}})\)
satisfies
\[
\|\overline{\mathbf{p}}-\overline{\mathbf{p}}^*\|
\le (1+c_\omega)\|\overline{\mathbf{m}}\| + (c_3+c_\omega c_2)\epsilon_1.
\]
Observe that
$
\mathbf{F}_{\rm map}(\overline{\mathbf{x}}^*,\overline{\mathbf{p}})
=\mathbf{G}_{\rm c}(\overline{\mathbf{p}}-\overline{\mathbf{p}}^*);
$
thus \(V_1(\overline{\mathbf{p}})\) is a positive
definite quadratic form and
\[
\|\nabla V_1(\overline{\mathbf{p}})\|
\le l_2\|\overline{\mathbf{p}}-\overline{\mathbf{p}}^*\|.
\]
Moreover, from the Lipschitz property of
\(\bm{\Phi}\) (constant \(l_1+l_\varepsilon/N\)),
\[
\|\mathbf{F}_{\rm map}(\overline{\mathbf{z}})
-\mathbf{F}_{\rm map}(\overline{\mathbf{x}}^*,\overline{\mathbf{p}})\|
\le l_3\|\overline{\mathbf{x}}-\overline{\mathbf{x}}^*\|.
\]
Insert
\[
\|\overline{\mathbf{x}}-\overline{\mathbf{x}}^*\|
\le c_0\|\overline{\mathbf{m}}\| + c_1\epsilon_1
\]
from Lemma~\ref{lem:main} and combine:
\[
\begin{aligned}
&\|\nabla V_1\|\,\|\mathbf{F}_{\rm map}(\overline{\mathbf{z}})
 -\mathbf{F}_{\rm map}(\overline{\mathbf{x}}^*,\overline{\mathbf{p}})\| \\
&\qquad\le
l_2\bigl[(1+c_\omega)\|\overline{\mathbf{m}}\|
 +(c_3+c_\omega c_2)\epsilon_1\bigr] \\
&\qquad\quad \times\,
l_3\bigl[c_0\|\overline{\mathbf{m}}\| + c_1\epsilon_1\bigr].
\end{aligned}
\]
Expanding the product and applying Young's inequality to the cross term yields
$
\|\nabla V_1\|\,\|\mathbf{F}_{\rm map}(\overline{\mathbf{z}})
-\mathbf{F}_{\rm map}(\overline{\mathbf{x}}^*,\overline{\mathbf{p}})\|
\le l_4\|\overline{\mathbf{m}}\|^2 + l_5\epsilon_1^2.
$
\end{proof}

We now assemble the auxiliary quantities for the
convergence statement. Define the Lyapunov
function
\begin{equation}
V(\overline{\mathbf{z}}) = V_1(\overline{\mathbf{p}}) + \frac{1}{2}\|\overline{\mathbf{x}}-\overline{\mathbf{x}}^*\|^2,
\label{eq:V_def}
\end{equation}
where $V_1$ is from
Corollary~\ref{cor:local_bounds}. Let
\begin{equation}
c_z = \frac{2(1+2c_0^2+2c_\omega^2)\,l_6}{\zeta} + 2\bigl(c_1^2+(c_3 + c_\omega c_2)^2\bigr),
\label{eq:cz_def}
\end{equation}
where $k_1 = l_2(1+c_\omega)^2 + c_0^2$, $k_2 =
l_2(c_3 + c_\omega c_2)^2 + c_1^2$, $\zeta =
(\eta-l_4)/\epsilon_1 -
(c_0+\gamma/2)/\epsilon_2$, $l_6 = l_5 +
c_1^2/(2\gamma\epsilon_2) + \zeta\,k_2/k_1$,
$\gamma >0 $.

\begin{theorem}
\label{thm:convergence} Suppose conditions
\textbf{(C1)}--\textbf{(C4)} hold. If $\zeta >
0$, then there exist constants
\(\epsilon_1^*>0\), \(c_1>0, c_2>0, c_3>0\) and a
finite time \(T \ge 0\) such that the inequality
\begin{equation}
\|\overline{\mathbf{z}}(t)-\overline{\mathbf{z}}^*\| \le \sqrt{c_z}\,\sqrt{\epsilon_1},
\label{eq:local_z_bound}
\end{equation}
holds for every $0<\epsilon_1<\epsilon_1^*$ and
$t\ge T$, where $c_z>0$ is defined
in~\eqref{eq:cz_def}. Moreover,
$\overline{\mathbf{x}}(t)$ \emph{converges
exponentially fast} to an
$O(\sqrt{\epsilon_1})$-neighborhood of
$\overline{\mathbf{x}}^*$ with rate
\(\zeta/(2k_1)\), satisfying
\begin{equation}
\|\overline{\mathbf{x}}(t)-\overline{\mathbf{x}}^*\| \le \sqrt{2V(T)}\,e^{-\frac{\zeta}{2k_1}(t-T)}
+ \sqrt{\frac{2k_1\,l_6}{\zeta}}\,\sqrt{\epsilon_1}.
\label{eq:local_x_bound}
\end{equation}
\end{theorem}

\begin{proof}
Since \(V_1\) is a positive definite quadratic
form in \(\overline{\mathbf{p}}-\overline{\mathbf{p}}^*\),
\(\|\nabla V_1(\overline{\mathbf{p}})\| \le
l_2\|\overline{\mathbf{p}}-\overline{\mathbf{p}}^*\|\).
Combining the pointwise bounds from
Lemmas~\ref{lem:main} and~\ref{lem:omega_bound}
yields
\[
V \le k_1\|\overline{\mathbf{m}}\|^2
+ k_2\,\epsilon_1^2,
\]
with $k_1$ and $k_2$ as defined
above.
The derivative of \(V\)
along~\eqref{eq:dynamics_global} admits the
decomposition
\[
\begin{aligned}
\frac{dV}{dt}
&=
\underbrace{(\nabla V_1)^{\mathsf T}
  \dot{\overline{\mathbf{p}}}\big|_{\overline{\mathbf{x}}
    =\overline{\mathbf{x}}^*}}_{\text{(i)}}
+ \frac{1}{\epsilon_2}
  \underbrace{(\overline{\mathbf{x}}-\overline{\mathbf{x}}^*)^{\mathsf T}
    \overline{\mathbf{m}}}_{\text{(iii)}} \\
&\qquad
+ \frac{1}{\epsilon_1}
  \underbrace{(\nabla V_1)^{\mathsf T}
    \bigl(\mathbf{F}_{\rm map}(\overline{\mathbf{z}})
      -\mathbf{F}_{\rm map}(\overline{\mathbf{x}}^*,
        \overline{\mathbf{p}})\bigr)}_{\text{(ii)}}.
\end{aligned}
\]
We estimate each term separately. From the
structure of \(\mathbf{G}_{\rm c}\),
\(\text{(i)}\le
-\frac{\eta}{\epsilon_1}\|\overline{\mathbf{m}}\|^2\).
By Corollary~\ref{cor:local_bounds},
\(\text{(ii)}\le l_4\|\overline{\mathbf{m}}\|^2 +
l_5\epsilon_1^2\). Invoking the bounds
\(\|\overline{\mathbf{x}}-\overline{\mathbf{x}}^*\|
\le c_0\|\overline{\mathbf{m}}\|+c_1\epsilon_1\)
(from Lemma~\ref{lem:main}) and
\(\|\overline{\bm{\omega}}-\overline{\bm{\omega}}^*\|
\le
c_\omega\|\overline{\mathbf{m}}\|+(c_3 + c_\omega c_2)\epsilon_1\),
and using Cauchy--Schwarz, we obtain
\(\text{(iii)}\le c_0\|\overline{\mathbf{m}}\|^2
+ c_1\epsilon_1\|\overline{\mathbf{m}}\|\).
Inserting these estimates and applying Young's
inequality
\[
c_1\epsilon_1\|\overline{\mathbf{m}}\|
\le \frac{\gamma}{2}\|\overline{\mathbf{m}}\|^2
+ \frac{c_1^2}{2\gamma}\epsilon_1^2
\]
with \(\gamma>0\) chosen as in the hypothesis, and
using \(\epsilon_1^2\le\epsilon_1\), we obtain
\begin{equation}
\frac{dV}{dt}
\le -\zeta\|\overline{\mathbf{m}}\|^2
+ \Bigl(l_5 + \frac{c_1^2}{2\gamma\epsilon_2}\Bigr)\epsilon_1,\qquad t\ge T.
\label{eq:Vdot_final}
\end{equation}

From~\eqref{eq:Vdot_final} and $l_6 \ge
l_5 + c_1^2/(2\gamma\epsilon_2)$, whenever
\(\|\overline{\mathbf{m}}\|^2\ge
2l_6\,\epsilon_1/\zeta\), we have
\(\dot{V}\le
-\frac{\zeta}{2}\|\overline{\mathbf{m}}\|^2\le0\).
Standard invariant set arguments then imply that
the trajectory eventually stays in
\[
\mathcal{S}
= \bigl\{\overline{\mathbf{z}}:
\|\overline{\mathbf{m}}\|^2
\le 2l_6\,\epsilon_1/\zeta\bigr\}.
\]
On \(\mathcal{S}\), using \((a+b)^2\le
2a^2+2b^2\) together with the pointwise bounds on
\(\overline{\mathbf{x}}\) and
\(\overline{\bm{\omega}}\) gives
\[
\begin{aligned}
\|\overline{\mathbf{z}}-\overline{\mathbf{z}}^*\|^2
&\le (1+2c_0^2+2c_\omega^2)\|\overline{\mathbf{m}}\|^2 \\
&\quad + 2\bigl(c_1^2+(c_3 + c_\omega c_2)^2\bigr)\epsilon_1^2
\le c_z\,\epsilon_1,
\end{aligned}
\]
with $c_z$ as defined
in~\eqref{eq:cz_def};
this proves~\eqref{eq:local_z_bound}.
Finally, inserting
\(\|\overline{\mathbf{m}}\|^2 \ge V/k_1
- (k_2/k_1)\epsilon_1^2\) into
\eqref{eq:Vdot_final} and using
\(\epsilon_1^2\le\epsilon_1\) yields
\[
\dot{V}
\le -\frac{\zeta}{k_1}V + l_6\,\epsilon_1.
\]
Gronwall's inequality then gives
\[
V(t)
\le V(T)e^{-\frac{\zeta}{k_1}(t-T)}
+ \frac{k_1\,l_6}{\zeta}\epsilon_1
\]
for all \(t\ge T\), and since \(V\ge
\frac12\|\overline{\mathbf{x}}-\overline{\mathbf{x}}^*\|^2\),
taking square roots
yields~\eqref{eq:local_x_bound}.
\end{proof}

By~\eqref{eq:local_z_bound}, the full state
\(\overline{\mathbf{z}}\) is ultimately bounded
within an \(O(\sqrt{\epsilon_1})\)-neighborhood
of the equilibrium. The \(\sqrt{\epsilon_1}\)
scaling reflects the quadratic structure of
\(V\), since the perturbation enters \(\dot{V}\)
at order \(O(\epsilon_1)\), the state-space
offset is \(O(\sqrt{\epsilon_1})\).
Equation~\eqref{eq:local_x_bound}
further decomposes the primal trajectory into a
decaying transient at exponential rate
\(\zeta/(2k_1)\) and a persistent steady-state
residual of order \(O(\sqrt{\epsilon_1})\).

\begin{corollary}
\label{cor:vanishing_limit}
Under the conditions of Theorem~\ref{thm:convergence},
as \(\epsilon_1 \to 0^+\) the ultimate bound
tends to zero, and the trajectory
exponentially approaches \(\overline{\mathbf{z}}^*\)
with rate \(\zeta/(2k_1)\).
In other words, the system recovers asymptotic
exponential stability in the vanishing-singular-perturbation
limit, with the quasi-steady-state manifold
coinciding exactly with the equilibrium manifold.
\end{corollary}

\begin{remark}
\label{rem:parameter_consolidation} According to
Theorem~\ref{thm:convergence}, the convergence
behavior of system~\eqref{eq:dynamics_global}
is fully characterized by two complementary
quantities, both explicitly parameterized by
the six tunable hyperparameters
\(\bm{\Theta} = (\epsilon_1, \epsilon_2, \eta,
\lambda, a, \varepsilon)\) together with the
problem data \((\overline{\mathbf{A}},
\overline{\mathbf{b}}, \mathbf{L})\):
\begin{itemize}
\item
  \textbf{Exponential convergence rate}
  \(\zeta/(2k_1)\), which controls how fast the
  transient decays. Larger \(\eta\) and
  \(\epsilon_2/\epsilon_1\) ratios accelerate
  convergence but must respect the stability
  condition \(\zeta > 0\), i.e.,
  \(\epsilon_2/\epsilon_1 > (c_0+\gamma/2)/(\eta-l_4)\)
  with \(\eta > l_4\).
\item
  \textbf{Ultimate bound} \(\sqrt{c_z\,\epsilon_1}\),
  which quantifies the steady-state precision.
  Reducing \(\epsilon_1\) tightens the bound at
  the cost of a smaller convergence rate,
  revealing a natural rate--accuracy trade-off
  intrinsic to the two-timescale architecture.
\end{itemize}
Although these expressions involve
numerous intermediate constants \(c_0, c_\omega,
l_1, \ldots, l_6, k_1, k_2, \zeta, c_z\), each is
explicitly determined by \(\bm{\Theta}\) and the
problem data. The feasible domain of \(\bm{\Theta}\) thus
defines a search space for global optimization,
which is explored by the bilevel PSO approach of
Section~\ref{sec:swarm}.
\end{remark}

\begin{remark}
\label{rem:two_phase_convergence}
Proposition~\ref{prop:boundary_layer} and
Theorem~\ref{thm:convergence} jointly
characterize the two-phase trajectory:
during \([0,T]\), \(\overline{\mathbf{p}}\)
converges exponentially to the slow manifold
at rate \(\alpha\) while
\(\overline{\mathbf{x}}\) varies by
\(O(\epsilon_1)\); during \([T,\infty)\),
\(\overline{\mathbf{x}}\) converges
exponentially to an
\(O(\sqrt{\epsilon_1})\)-neighborhood of
\(\overline{\mathbf{x}}^*\) at rate
\(\zeta/(2k_1)\). Both phases are verified
experimentally in
Section~\ref{sec:simulation}.
\end{remark}

\section{Bilevel PSO-Neurodynamic Optimization
with RLGL Acceleration}
\label{sec:swarm}

The collaborative neurodynamic multi-agent
system~\eqref{eq:dynamics_global} guarantees only
local convergence within the basin prescribed by
Corollary~\ref{cor:gradient_growth} and
condition \textbf{(C4)}; outside this
region, convergence to spurious local minima may
occur.  This limitation has motivated a bilevel
approach, extending the neurodynamic bilevel
optimization paradigm~\cite{qin2017neurodynamic}
to the nonconvex sparse recovery setting, that
searches over a six-dimensional
hyperparameter vector \(\bm{\Theta} =
(\epsilon_1, \epsilon_2, \eta, \lambda, a,
\varepsilon) \in \mathcal{D} \subset
\mathbb{R}^6\), where \(\mathcal{D}\) is a
compact hyper-rectangle encoding physically
admissible bounds.  Specifically,
\(\epsilon_1\) is confined to
\((0,\epsilon_1^*]\) so that the
singular-perturbation threshold of
Theorem~\ref{thm:convergence} is respected;
\(\epsilon_2\) and \(\eta\) occupy the region
where the stability
condition~\eqref{eq:stability_check} holds
with a prescribed margin; \(\lambda\) spans a
range bracketing the noise level of
\(\mathbf{b}_i\); \(a\) is restricted to
\((2,\infty)\), the interval on which SCAD is
concave in the penalty; and \(\varepsilon\) is
bounded above so that the mollified penalty
retains the SCAD transition structure
(Section~\ref{sec:problem}).
Within this approach, the
inner layer integrates
system~\eqref{eq:dynamics_global} to stationarity
under fixed \(\bm{\Theta}\), while the outer
layer performs derivative-free search to optimize
recovery performance subject to the stability
constraints derived below.

The extended Bayesian information criterion
(eBIC)~\cite{chen2008extended} is adopted as the
outer objective
\begin{equation}
  \label{eq:ebic}
  \mathcal{J}(\bm{\Theta}) = m \log(\operatorname{RSS}/m) + \hat{s} \log m + 2\gamma \log \binom{n}{\hat{s}},
\end{equation}
where \(\operatorname{RSS} =
\sum_{i=1}^{N}\|\mathbf{A}_i
\mathbf{x}_i(T)-\mathbf{b}_i\|^2\),
\(\bar{\mathbf{x}}(T)\) denotes the consensus
average at terminal time \(T\), and \(\hat{s} =
\|\bar{\mathbf{x}}(T)\|_{0,\delta}\) is the
\(\delta\)-truncated sparsity level, i.e., the
number of entries of \(\bar{\mathbf{x}}(T)\)
whose magnitude exceeds a small threshold
\(\delta>0\).  The
parameter \(\gamma \in [0,1]\) controls
complexity penalization.  Unlike the ordinary BIC,
the combinatorial term
\(2\gamma \log \binom{n}{\hat{s}}\) penalizes the
number of candidate models of size \(\hat{s}\),
which is essential in the high-dimensional
regime \(n\gg m\) where BIC tends to overselect
spurious variables.

The bilevel structure follows directly from the
convergence analysis.
Theorem~\ref{thm:convergence} guarantees
semiglobal practical convergence to an
\(O(\sqrt{\epsilon_1})\) neighborhood of a
stationary point, yet two fundamental issues
persist. The nonconvex SCAD landscape admits
multiple equilibria with different support
sets~\cite{fan2001variable}, and the collaborative neurodynamic multi-agent system
may converge to a poor local minimum depending on
\((\lambda,a,\varepsilon)\); moreover, the
residual radius cannot be made arbitrarily small
without incurring prohibitive
stiffness~\cite{khalil2002nonlinear}.  Hence, the
collaborative neurodynamic multi-agent system alone only provides local
convergence for a fixed configuration, and global
parameter selection remains necessary. The outer
PSO-RLGL optimization layer addresses this need by exploring the
discontinuous mapping \(\bm{\Theta} \mapsto
\hat{s}\), where small changes in \(\lambda\) or
\(\varepsilon\) can abruptly alter the active
set.  The eBIC criterion serves as a selection
mechanism among the resulting solution manifolds,
and the derivative-free, population-based nature of
PSO enables traversal of this jump-discontinuous landscape
without gradient information~\cite{clerc2002particle}.

Importantly, the convergence theory does not
merely coexist with the PSO; it actively
constrains the search.  A sufficient condition
for the parameter requirement in
Theorem~\ref{thm:convergence} is given by
\begin{equation}
\label{eq:stability_check}
\eta > \frac{\|\mathbf{L}\otimes\mathbf{I}_n\|}{\epsilon_2} +  l_1 + l_\varepsilon/N,
\end{equation}
where \(l_1 = \max_i \|\mathbf{A}_i^{\mathsf{T}}\mathbf{A}_i\|\) is the spectral-norm bound on the local Gram matrices.  The three right-hand-side terms bound the consensus coupling, the least-squares gradient, and the per-agent smoothed SCAD gradient, respectively; \eqref{eq:stability_check} thus enforces that the descent gain \(\eta\) dominates all destabilizing contributions.  Particles violating~\eqref{eq:stability_check}
are rejected, confining the computational budget
to provably stable regions.  Moreover,
Theorem~\ref{thm:convergence} explicitly links
the residual radius to \(\epsilon_1\) through
\(\zeta\) and \(l_6\); consequently,
\(\epsilon_1\) becomes a decision variable
directly controlling solution accuracy, and the
PSO searches for the optimal trade-off between
convergence speed and precision.

The analysis further reveals a functional
dichotomy among the parameters.  The time-scale
parameters \((\epsilon_1,\epsilon_2,\eta)\) govern
convergence rate and precision, but once
\eqref{eq:stability_check} is satisfied with
sufficient margin, further tuning yields
diminishing returns.  In contrast, the SCAD
parameters \((\lambda,a)\) and the mollification
parameter \(\varepsilon\) directly reshape the
nonconvex penalty landscape, inducing persistent,
nonlinear, and often multimodal responses in
\(\mathcal{J}\).  The RLGL mechanism~\cite{avrachenkov2022redlight} exploits
this asymmetry by dynamically identifying
dimensions whose sensitivity to \(\mathcal{J}\)
has saturated and concentrating computational
resources on those that remain influential. To
this end, the first-order sensitivity is estimated via
central finite differences
\begin{equation}
\label{eq:sensitivity}
s_{p,d} = \big| \mathcal{J}(\bm{\Theta}_p + \delta_d \mathbf{e}_d) - \mathcal{J}(\bm{\Theta}_p - \delta_d \mathbf{e}_d) \big|,
\end{equation}
with a dimension-adaptive perturbation \(\delta_d>0\).
The perturbation is scaled to the search range of
dimension \(d\), \(\delta_d = \kappa
|\mathcal{D}_d|\) with \(\kappa\ll 1\), ensuring
that finite-difference estimates are comparable
across parameters whose nominal scales differ by
orders of magnitude.  To reduce overhead, sensitivity estimation is
performed only for a sampled subset of \(q\cdot
N_p\) particles per generation, \(q\in(0,1]\).
The green-light set for particle \(p\) is then
determined by a relative threshold
\begin{equation*}
\begin{aligned}
\mathcal{G}_p &= \bigl\{ d \in \{1,\ldots,6\} : s_{p,d} \ge \tau^{(k)} s_{p,\max} \bigr\},\\
s_{p,\max} &= \max_{d'} s_{p,d'}.
\end{aligned}
\end{equation*}
Only dimensions in \(\mathcal{G}_p\) are updated;
the velocities of the remaining dimensions are
frozen, effectively skipping the corresponding
evaluations until reactivation.  A global
all-green regeneration cycle is enforced every
\(K_{\text{reset}}\) iterations to preserve the
ergodic exploration properties of the baseline
PSO~\cite{clerc2002particle} and prevent
permanently frozen dimensions from trapping the
swarm in a low-dimensional subspace.  The threshold
follows a linear annealing schedule \(\tau^{(k)}
= \tau_{\max} - \frac{k}{K_{\max}}(\tau_{\max} -
\tau_{\min})\), reflecting the structural insight
that early search should focus on the most
sensitive (typically SCAD) parameters, while the
decaying threshold gradually reactivates time-scale
parameters for fine-tuning.

The expected per-generation cost of RLGL, with average
active ratio \(\bar{\rho}\) and sampling rate
\(q\), reduces from the baseline
\(N_p(6C_{\text{sens}}+C_{\text{int}})\) to
\(N_p(q\cdot6\bar{\rho}C_{\text{sens}}+C_{\text{int}})\),
where \(C_{\text{int}}\) is the cost of one inner
integration and
\(C_{\text{sens}}=2C_{\text{int}}\).  This yields
a significant computational saving without
compromising the convergence behavior established
in Theorem~\ref{thm:convergence}.

The outer search is implemented via the standard
constriction-factor PSO~\cite{clerc2002particle}.
PSO is preferred over gradient-based alternatives
because \(\bm{\Theta}\mapsto\mathcal{J}(\bm{\Theta})\)
is nonsmooth and jump-discontinuous at support
transitions, making gradient estimators unreliable;
its population-based structure simultaneously
maintains candidates across multiple
basins~\cite{kennedy1995pso, clerc2002particle}.
Each particle \(p\) maintains its position,
velocity, personal best \(\mathbf{p}_p\), and the
swarm tracks the global best \(\mathbf{g}\).  The
update rule is
\begin{align}
\mathbf{v}_p &\leftarrow \chi \big[ \mathbf{v}_p + c_1 r_1 (\mathbf{p}_p - \bm{\Theta}_p) + c_2 r_2 (\mathbf{g} - \bm{\Theta}_p) \big], \nonumber\\
\bm{\Theta}_p &\leftarrow \Pi_{\mathcal{D}}(\bm{\Theta}_p + \mathbf{v}_p),
\label{eq:pso}
\end{align}
with \(\chi = 0.729\), \(c_1 = c_2 = 2.05\),
\(r_1,r_2\sim\mathcal{U}(0,1)\), and
\(\Pi_{\mathcal{D}}\) the Euclidean projection
onto \(\mathcal{D}\).
Algorithm~\ref{alg:bilevel} summarizes the
complete procedure.

\begin{algorithm}[htbp]
\caption{RLGL-Accelerated Bilevel PSO-Neurodynamic Optimization}
\label{alg:bilevel}
\begin{algorithmic}[1]
\REQUIRE Data $(\mathbf{A}_i,\mathbf{b}_i)_{i=1}^N$, graph $\mathcal{G}$, swarm size $N_p$, max generations $K_{\max}$, reset period $K_{\text{reset}}$, horizon $T$, thresholds $\tau_{\max},\tau_{\min}$, sampling rate $q$
\ENSURE Optimal $\bm{\Theta}^*$ and sparse solution $\mathbf{x}^*$
\STATE Initialize $\bm{\Theta}_p\sim\mathcal{U}(\mathcal{D})$, $\mathbf{v}_p\gets\mathbf{0}$, $p=1,\dots,N_p$
\FOR{$p=1$ to $N_p$}
    \IF{\eqref{eq:stability_check} is violated}
        \STATE $\mathcal{J}(\bm{\Theta}_p)\gets\infty$
    \ELSE
        \STATE Integrate~\eqref{eq:dynamics_global}, evaluate $\mathcal{J}(\bm{\Theta}_p)$ via~\eqref{eq:ebic}
    \ENDIF
    \STATE $\mathbf{p}_p\gets\bm{\Theta}_p$
\ENDFOR
\STATE $\mathbf{g}\gets\arg\min_p \mathcal{J}(\bm{\Theta}_p)$
\FOR{$k=1$ to $K_{\max}$}
    \STATE $\tau^{(k)}\gets\tau_{\max} - \frac{k}{K_{\max}}(\tau_{\max}-\tau_{\min})$
    \FOR{$p=1$ to $N_p$}
        \IF{$\mathrm{rand}()<q$ \AND $k\bmod K_{\text{reset}}\neq0$}
            \STATE Compute $s_{p,d}$ via~\eqref{eq:sensitivity}, $\mathcal{G}_p\gets\{d:s_{p,d}\ge\tau^{(k)}s_{p,\max}\}$
        \ELSE
            \STATE $\mathcal{G}_p\gets\{1,\dots,6\}$
        \ENDIF
        \STATE Update $\bm{\Theta}_p$ using~\eqref{eq:pso} only for $d\in\mathcal{G}_p$, freeze other velocities
        \IF{$\bm{\Theta}_p$ changed}
            \IF{\eqref{eq:stability_check} is violated}
                \STATE $\mathcal{J}(\bm{\Theta}_p)\gets\infty$
            \ELSE
                \STATE Integrate~\eqref{eq:dynamics_global}, evaluate $\mathcal{J}(\bm{\Theta}_p)$
            \ENDIF
        \ENDIF
        \STATE Update $\mathbf{p}_p$ and $\mathbf{g}$
    \ENDFOR
\ENDFOR
\STATE Refine: integrate with $\bm{\Theta}^*=\mathbf{g}$ to convergence; $\mathbf{x}^*\gets\frac1N\sum_i\mathbf{x}_i(T)$
\RETURN $\bm{\Theta}^*$, $\mathbf{x}^*$
\end{algorithmic}
\end{algorithm}

\begin{remark}
\label{rem:sparsity_pso}
The structure-preserving mollification of
Section~\ref{sec:problem}
analytically preserves the
unbiasedness component of the SCAD
oracle property
(Remark~\ref{rem:unbiasedness_preservation}),
but the sparsity selection
component---which requires exact support
recovery in the noiseless limit---is
not guaranteed by mollification alone,
since the smoothed penalty remains
strictly positive off the support.
Theorem~\ref{thm:convergence} clarifies
that the inner system converges only to
an \(O(\sqrt{\epsilon_1})\)-neighborhood
of a stationary point, leaving a residual
radius that may prevent exact thresholding.
The PSO--RLGL outer layer closes this gap
empirically: by treating
\((\lambda,a,\varepsilon)\) as decision
variables and selecting configurations via
eBIC~\eqref{eq:ebic}, the bilevel search
identifies parameter settings under which
the terminal solution attains the correct
support.  The feasibility constraint
\(\epsilon_1<\epsilon_1^*\)
(Theorem~\ref{thm:convergence}) guarantees
that the ultimate bound
\(\sqrt{c_z\epsilon_1}\) is well defined,
and the PSO search keeps \(\epsilon_1\)
sufficiently small so that this bound
stays below the smallest nonzero
coefficient magnitude, ensuring that the
post-processing hard-thresholding step
yields exact support recovery.  The empirical
validation in
Section~\ref{sec:simulation} confirms
this prediction across all four
scenarios.
\end{remark}

\begin{remark}
\label{rem:global_convergence}
The bilevel framework converges in
probability to the global optimum of
\(\mathcal{J}\) over \(\mathcal{D}\).
Since \(\mathcal{D}\) is compact and
\(\mathcal{J}\) is bounded below, the
global infimum \(\mathcal{J}^*\) exists.
Theorem~\ref{thm:convergence} ensures
that \(\mathcal{J}(\bm{\Theta})\) is
well-defined for every feasible
\(\bm{\Theta}\).
The stochastic nature of the PSO
update~\eqref{eq:pso}, with i.i.d.\
random coefficients
\(r_1,r_2\sim\mathcal{U}(0,1)\), together
with the periodic all-green reset,
provides a mechanism for sustained
exploration of \(\mathcal{D}\).
Under standard exploration
assumptions for stochastic search on
compact
sets~\cite{clerc2002particle}, one
obtains, for any
\(\varepsilon_g>0\),
\(\lim_{K\to\infty}
\Pr[\mathcal{J}(\mathbf{g}^{(K)})
- \mathcal{J}^* < \varepsilon_g] = 1\).
This parallels the probabilistic global
convergence result
of~\cite{xia2024momentum}.
\end{remark}

\section{Numerical Examples}
\label{sec:simulation}
\begingroup\sloppy\emergencystretch=3em
We evaluate all methods on four synthetic scenarios
of increasing difficulty: Low-D
($n{=}50$, $m{=}40$, $s{=}5$, $\sigma{=}0.01$) as an
easy recovery regime, Mid-D
($n{=}150$, $m{=}120$, $s{=}20$, $\sigma{=}0.02$) as the
nominal operating point, High-D
($n{=}300$, $m{=}240$, $s{=}50$, $\sigma{=}0.05$)
stressing high-dimensional recovery, and
High-Noise
($n{=}150$, $m{=}120$, $s{=}20$, $\sigma{=}0.1$) probing
noise robustness.  In each scenario the
ground-truth signal
\(\mathbf{x}_0\in\mathbb{R}^n\) contains \(s\)
nonzero entries drawn i.i.d.\ from
\(\mathcal{N}(0,1)\) with support
\(\mathcal{S}_0\) uniformly random over
\(\{1,\dots,n\}\); the sensing matrix
\(\mathbf{A}\in\mathbb{R}^{m\times n}\) has
i.i.d.\ Gaussian entries scaled by
\(1/\sqrt{m}\); and the measurements
\(\mathbf{b}=\mathbf{A}\mathbf{x}_0+\boldsymbol{\eta}\)
are corrupted by Gaussian noise
\(\boldsymbol{\eta}\sim\mathcal{N}(0,\sigma^2\mathbf{I})\).
All distributed methods operate on a 5-agent ring
network with even measurement partition.  To
ensure that performance differences reflect only
the algorithmic paradigm, all methods share the
same post-processing hard-thresholding step
\(\max\{0.1\hat{\lambda},\,0.001\max_j|\hat{x}_j|\}\)
and the same composite fitness
\(\mathrm{Rel.Err}+1.5|\hat{s}-s|/s+(1-\mathrm{precision})\).
The SCAD shape parameter \(a\in[2.1,3.9]\) is
jointly tuned by PSO for the neurodynamic methods
and fixed at \(3.0\) for the decentralized
smoothing ADMM (DSAD).
\par\endgroup

\subsection{Baseline Comparison}
\label{sec:neurodynamic_comparison}

To assess the proposed approach, we compare with
three distributed neurodynamic
baselines---the augmented-Lagrangian
projection neural network (PNN-AL)
\cite{xia2023collaborative}, the
two-timescale projection neural network
(Two-TS PNN)
\cite{huang2024twotimescale}, and the
penalty neurodynamic approach (Penalty-ND)
\cite{jia2024collective}---and one iterative
reference, DSAD \cite{mirzaeifard2025decentralized},
which solves the \emph{same SCAD-regularized
problem} by operator splitting rather than by
neurodynamic evolution.  The baselines are tuned
by an identical PSO configuration (swarm size 20,
20 iterations) to ensure fairness, while DSAD is
tuned by grid search over its regularization
parameters.

Table~\ref{tab:neurodynamic_comparison} reports
the relative error
(\(\|\hat{\mathbf{x}}-\mathbf{x}_0\|_2/\|\mathbf{x}_0\|_2\)),
estimated sparsity
(\(|\hat{\mathcal{S}}|\)), precision
(\(|\hat{\mathcal{S}}\cap\mathcal{S}_0|/|\hat{\mathcal{S}}|\)),
and recovery rate
(\(|\hat{\mathcal{S}}\cap\mathcal{S}_0|/|\mathcal{S}_0|\))
averaged over 20 independent realizations per
scenario.  The proposed method is the
\emph{best-performing neurodynamic approach}
across all scenarios, and three findings
directly map to the theoretical contributions.
\emph{(i) Sparsity recovery (Contribution 1):}
the proposed method is the \emph{only}
neurodynamic approach whose estimated
sparsity exactly matches the true sparsity in
all four scenarios, while the projection-based
methods (PNN-AL, Two-TS PNN)
consistently underestimate sparsity in the
Low-D and Mid-D regimes.  As anticipated by
Remark~\ref{rem:sparsity_pso}, sparsity
recovery is achieved empirically through the
PSO--RLGL outer search, with the mollification
parameter \(\varepsilon\) constrained by the
local-optimality condition of
Theorem~\ref{thm:local_optimality} and the
timescale parameter \(\epsilon_1\) constrained
by \(\epsilon_1^*\) to keep the residual below
the smallest nonzero coefficient, while the
unbiasedness on large coefficients is guaranteed
analytically by
Remark~\ref{rem:unbiasedness_preservation}.
\emph{(ii) Accuracy (Contributions 1--2):}
in Low-D and Mid-D the proposed method achieves
perfect precision and recovery (\(1.00\)) with
the lowest relative error, matching the
iterative DSAD reference and surpassing all
competing neurodynamic methods by 5--22 times
in accuracy.  \emph{(iii) Theory-guided precision
(Theorem~\ref{thm:local_optimality}):} in the
High-D regime, reducing the
mollification parameter \(\varepsilon\) to
\([10^{-5},0.1]\), consistent with the
local-optimality condition, raises precision
from \(0.62\) to \(0.94\) while keeping the
sparsity exactly at the true value \(50\).

\begin{table}[t]
\centering
\caption{Comparison on four synthetic scenarios.
Best per scenario in bold.}
\label{tab:neurodynamic_comparison}
\small
\begin{tabular}{@{}lcccc@{}}
\hline
\textbf{Method} & \textbf{Rel.Err} & \textbf{Sparsity} & \textbf{Prec.} & \textbf{Rec.} \\
\hline
\multicolumn{5}{c}{\textit{Low-D ($n{=}50$, $s{=}5$)}} \\
\textbf{Proposed}            & \textbf{0.005} & \textbf{5} & \textbf{1.00} & \textbf{1.00} \\
PNN-AL                                     & 0.114 & 3  & 1.00 & 0.60 \\
Two-TS PNN                                & 0.114 & 3  & 1.00 & 0.60 \\
Penalty-ND                                & 0.043 & 5  & 1.00 & 1.00 \\
DSAD                                      & 0.005 & 5  & 1.00 & 1.00 \\
\hline
\multicolumn{5}{c}{\textit{Mid-D ($n{=}150$, $s{=}20$)}} \\
\textbf{Proposed}            & \textbf{0.017} & \textbf{20} & \textbf{1.00} & \textbf{1.00} \\
PNN-AL        & 0.089 & 17 & 1.00 & 0.85 \\
Two-TS PNN   & 0.090 & 17 & 1.00 & 0.85 \\
Penalty-ND       & 0.255 & 17 & 1.00 & 0.85 \\
DSAD  & 0.017 & 20 & 1.00 & 1.00 \\
\hline
\multicolumn{5}{c}{\textit{High-D ($n{=}300$, $s{=}50$)}} \\
\textbf{Proposed}            & 0.069 & \textbf{50} & 0.94 & 0.94 \\
PNN-AL                                     & 0.192 & 32 & \textbf{1.00} & 0.64 \\
Two-TS PNN                                & 0.242 & 40 & 0.85 & 0.68 \\
Penalty-ND                                & 0.112 & 53 & 0.89 & 0.94 \\
DSAD                                      & \textbf{0.065} & 49 & 0.96 & 0.94 \\
\hline
\multicolumn{5}{c}{\textit{High-Noise ($n{=}150$, $s{=}20$, $\sigma{=}0.1$)}} \\
\textbf{Proposed}            & 0.168 & \textbf{20} & 0.90 & 0.90 \\
PNN-AL                                     & \textbf{0.163} & 21 & 0.90 & 0.95 \\
Two-TS PNN                                & 0.168 & 19 & \textbf{0.95} & 0.90 \\
Penalty-ND                                & 0.721 & 66 & 0.29 & 0.95 \\
DSAD                                      & 0.176 & 21 & 0.90 & 0.95 \\
\hline
\end{tabular}
\end{table}

\begin{figure}[t]
\centering
\includegraphics[width=\columnwidth]{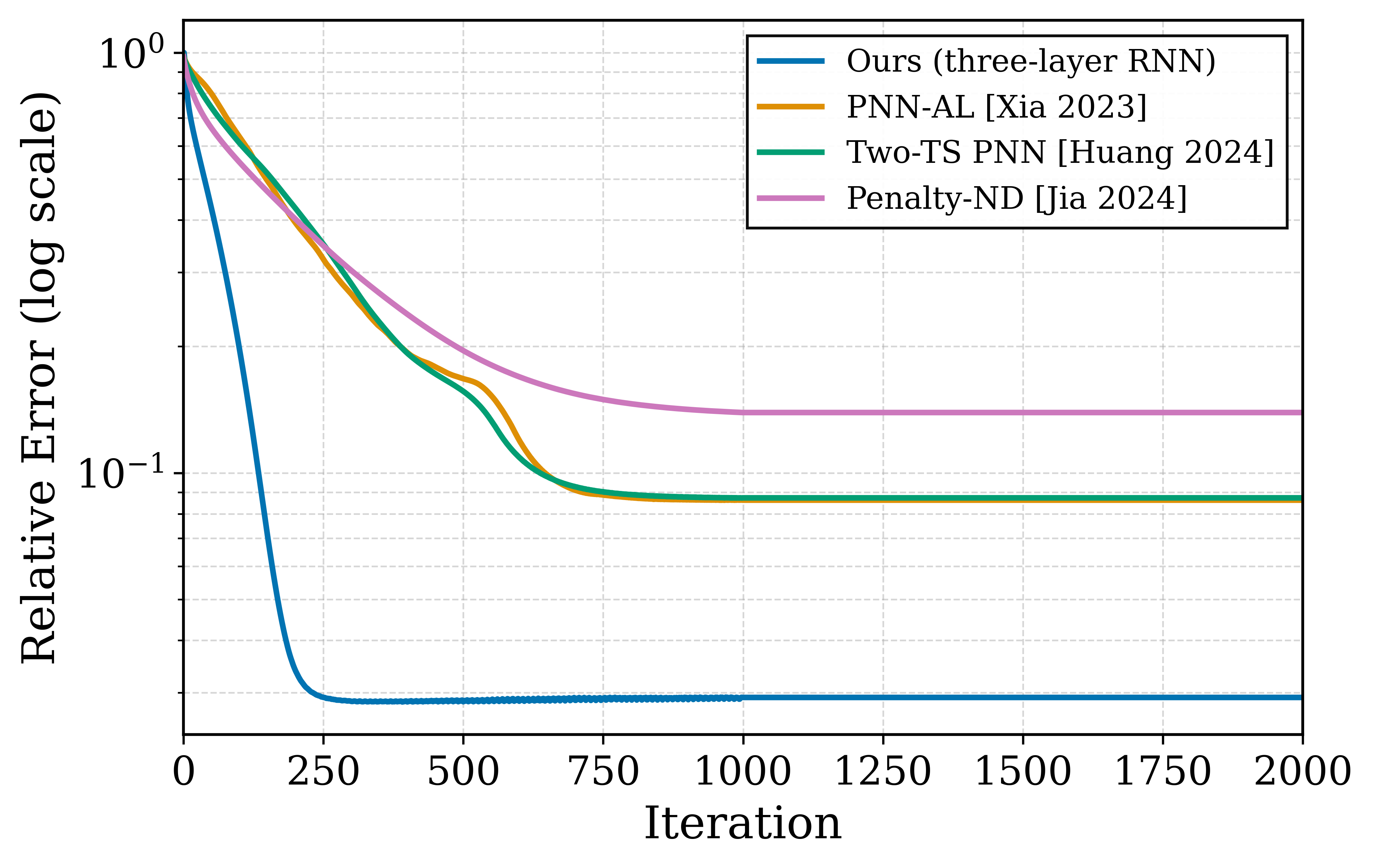}
\caption{Neurodynamic convergence on the Mid-D scenario.}
\label{fig:neurodynamic_convergence}
\end{figure}

The convergence curves in
Figure~\ref{fig:neurodynamic_convergence} on the
Mid-D scenario further highlight the architectural
advantage.  Because DSAD's per-iteration error
already incorporates the SCAD proximal step, a
structurally different convergence measure from
the recurrent neural network (RNN) state
evolution of neurodynamic methods, the
comparison is restricted to the four neurodynamic
approaches.  The proposed three-layer RNN with
leaky-integrator momentum converges roughly three
times faster than PNN-AL, Two-TS PNN, and
Penalty-ND, and reaches a steady-state error
floor of \(\sim 3{\times}10^{-2}\), 3--5 times
lower than the competing neurodynamic methods,
consistent with the exponential
ultimate-boundedness rate of
Theorem~\ref{thm:convergence}.

\subsection{RLGL Ablation}
\label{sec:rlgl_ablation}

Unlike standard PSO, which treats all
hyperparameters homogeneously, RLGL exploits the
\emph{functional dichotomy} of the six system
parameters.  The timescale parameters
\((\epsilon_1, \epsilon_2, \eta)\) govern
convergence speed, while the landscape-shaping
parameters \((\lambda, a, \varepsilon)\) govern
support selection.  RLGL freezes timescale
dimensions whose sensitivity falls below a
threshold, concentrating the swarm budget on the
active landscape parameters.

To isolate the contribution of each design
choice, we compare four PSO variants on the
Mid-D scenario with identical, reduced swarm
size (15 particles) and iteration budget (12
iterations)---smaller than the main comparison
to stress-test the freezing mechanism under
tight budgets---averaged over three random
seeds:
\begin{itemize}
\item \textbf{RLGL-PSO (ours):} functional
    dichotomy + sensitivity-based freezing;
\item \textbf{Standard PSO:} no freezing, all
    dimensions updated every iteration;
\item \textbf{Random Freeze PSO:} dichotomy
    preserved but freezing decisions are random
    (ablates the sensitivity analysis);
\item \textbf{All-Freeze PSO:} sensitivity-based
    freezing applied to all six dimensions
    (ablates the functional dichotomy).
\end{itemize}

\begin{figure}[t]
\centering
\includegraphics[width=\columnwidth]{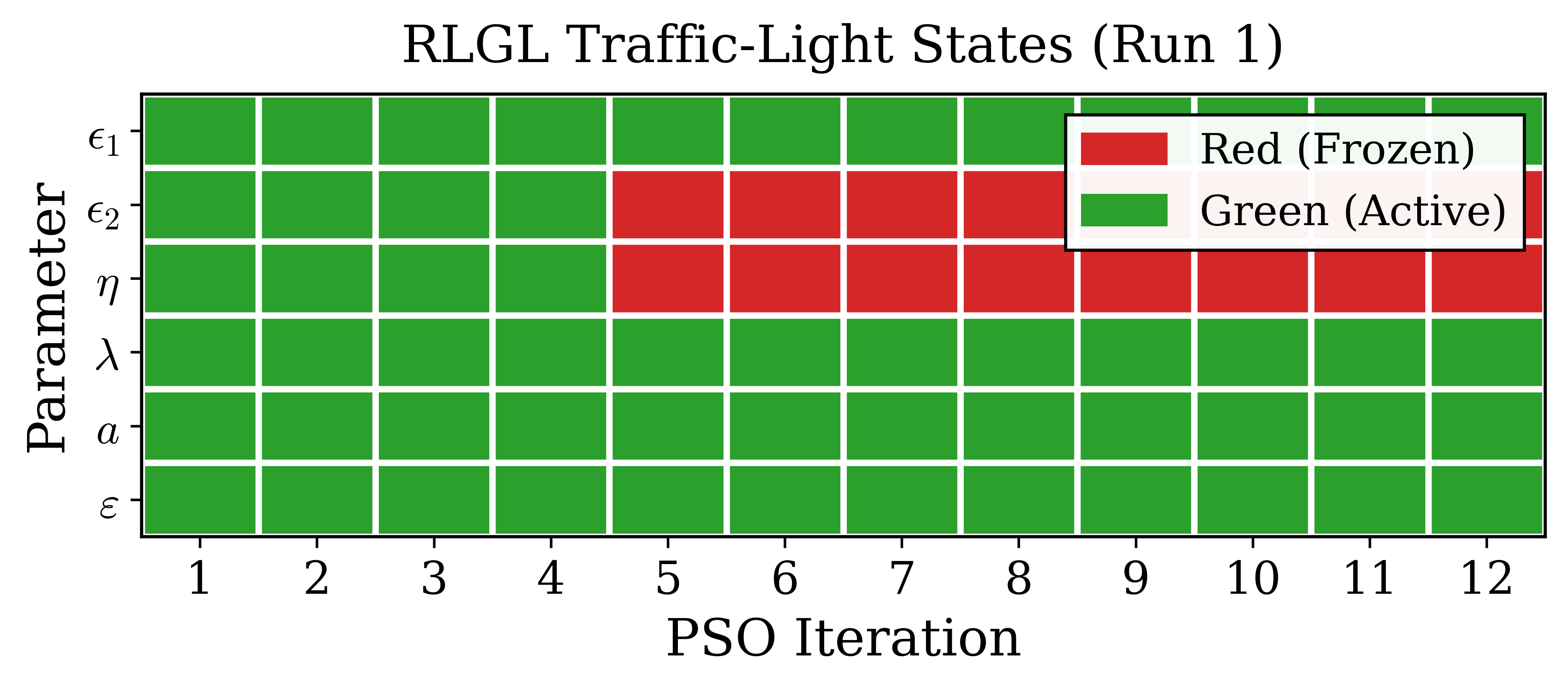}
\caption{RLGL traffic-light states across PSO iterations.}
\label{fig:rlgl_states}
\end{figure}

Figure~\ref{fig:rlgl_states} visualizes the
traffic-light states of RLGL-PSO across iterations.
The timescale parameters \(\epsilon_2\) and
\(\eta\) are frozen (red) from iteration~5 onward
once their sensitivity saturates, while the
landscape parameters \(\lambda, a, \varepsilon\)
remain green throughout.  This selective freezing
reduces the effective search dimensionality by
one-third without degrading recovery quality.

\begin{figure}[t]
\centering
\includegraphics[width=\columnwidth]{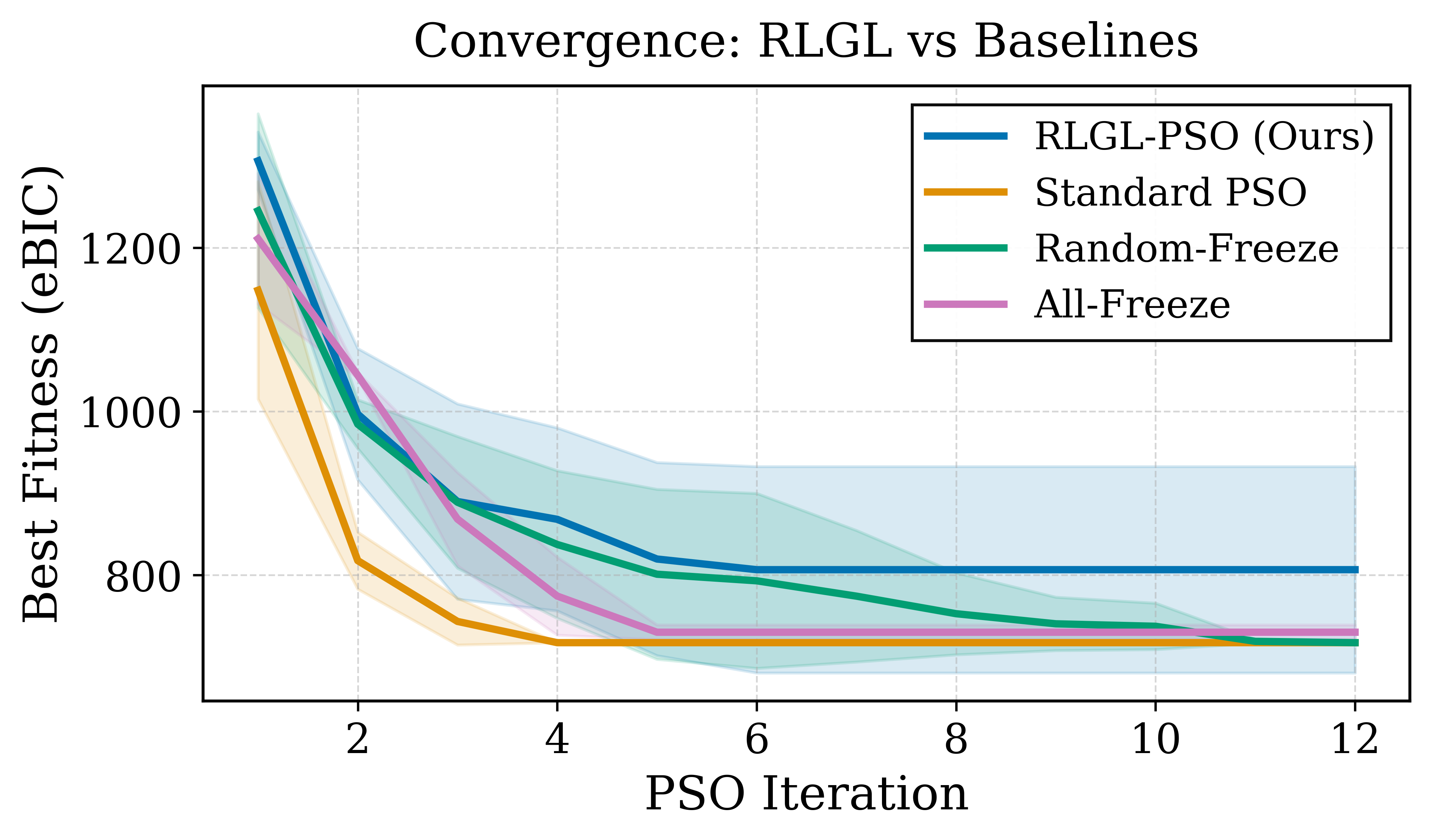}
\caption{Convergence of four PSO variants.}
\label{fig:ablation_conv}
\end{figure}

The convergence curves in
Figure~\ref{fig:ablation_conv} show that RLGL-PSO
achieves competitive fitness with fewer effective
dimension updates.

\begin{figure}[t]
\centering
\includegraphics[width=\columnwidth]{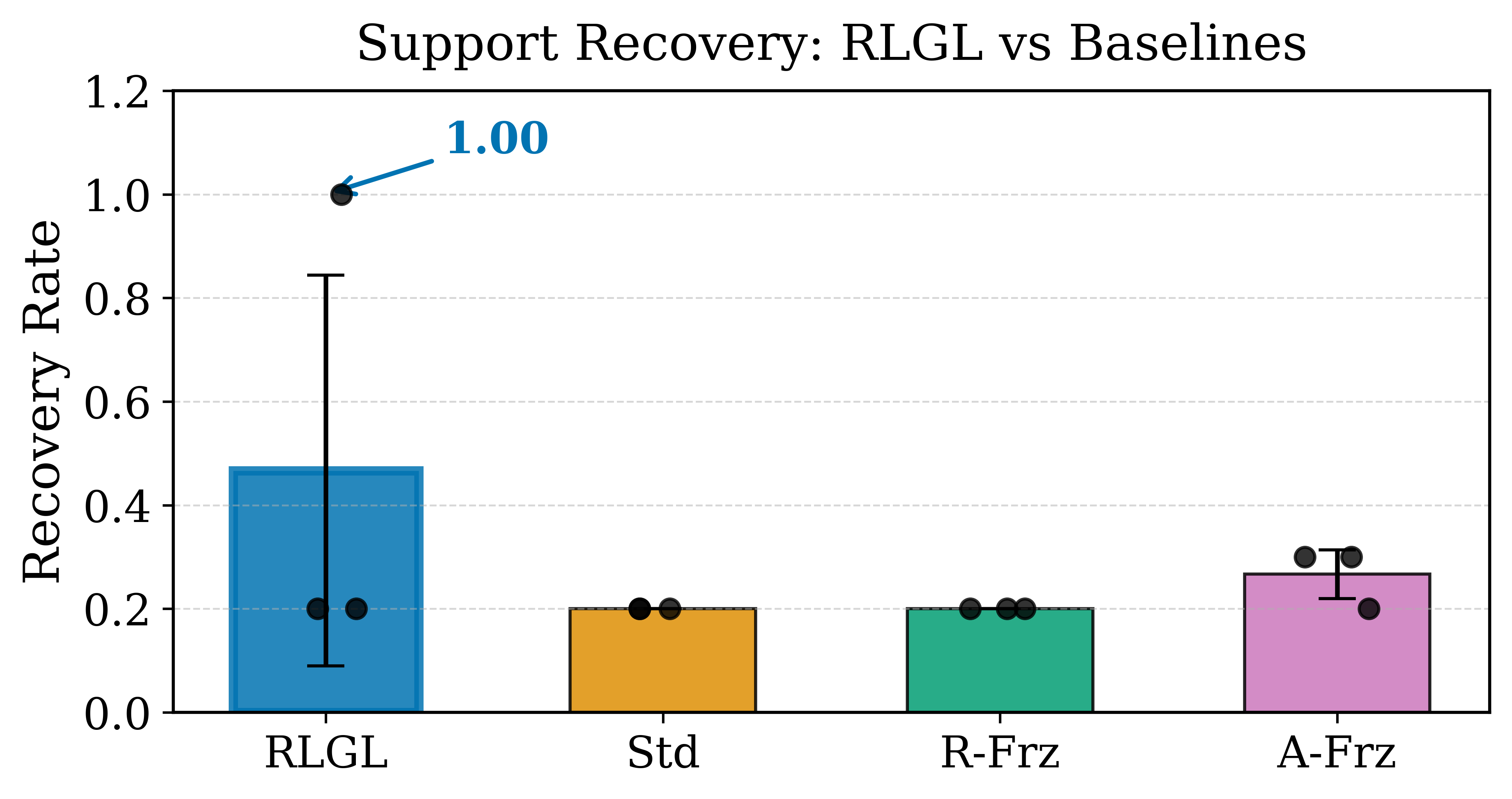}
\caption{Support recovery rate across four PSO variants.}
\label{fig:rlgl_ablation_bar}
\end{figure}

Figure~\ref{fig:rlgl_ablation_bar} quantifies the
RLGL contribution across the four variants.
RLGL-PSO is the only variant achieving perfect
support recovery (recovery rate \(=1.00\)) in any
of the three runs, whereas Standard PSO and
Random-Freeze PSO remain at \(0.20\) across all
runs, and All-Freeze PSO, which drops the
functional dichotomy, reaches only \(0.27\).
This confirms that the \emph{combination} of
sensitivity-based freezing and functional
dichotomy is essential.

\subsection{Convergence rate verification}
\label{sec:theorem_validation}

Theorem~\ref{thm:convergence} and
Proposition~\ref{prop:boundary_layer} predict
two-phase exponential convergence with
explicitly parameterized behavior: a slow-phase
rate \(\zeta/(2k_1)\) that is
\(\epsilon_1\)-independent, a boundary layer
duration \(T \approx 5\epsilon_1/\alpha\), and an
ultimate bound \(\sqrt{c_z\epsilon_1}\) that
increases with \(\epsilon_1\).  We verify these
predictions on a synthetic sparse recovery problem
(\(n=100, s=10, \sigma=0.01\)) by fixing all
parameters at their RLGL-PSO--tuned optimal values
and sweeping \(\epsilon_1 \in \{0.1, 0.3, 0.5,
1.0, 2.0, 5.0, 10.0\}\epsilon_1^*\), where
\(\epsilon_1^*\) is the optimal value.

\begin{figure}[t]
\centering
\includegraphics[width=\columnwidth]{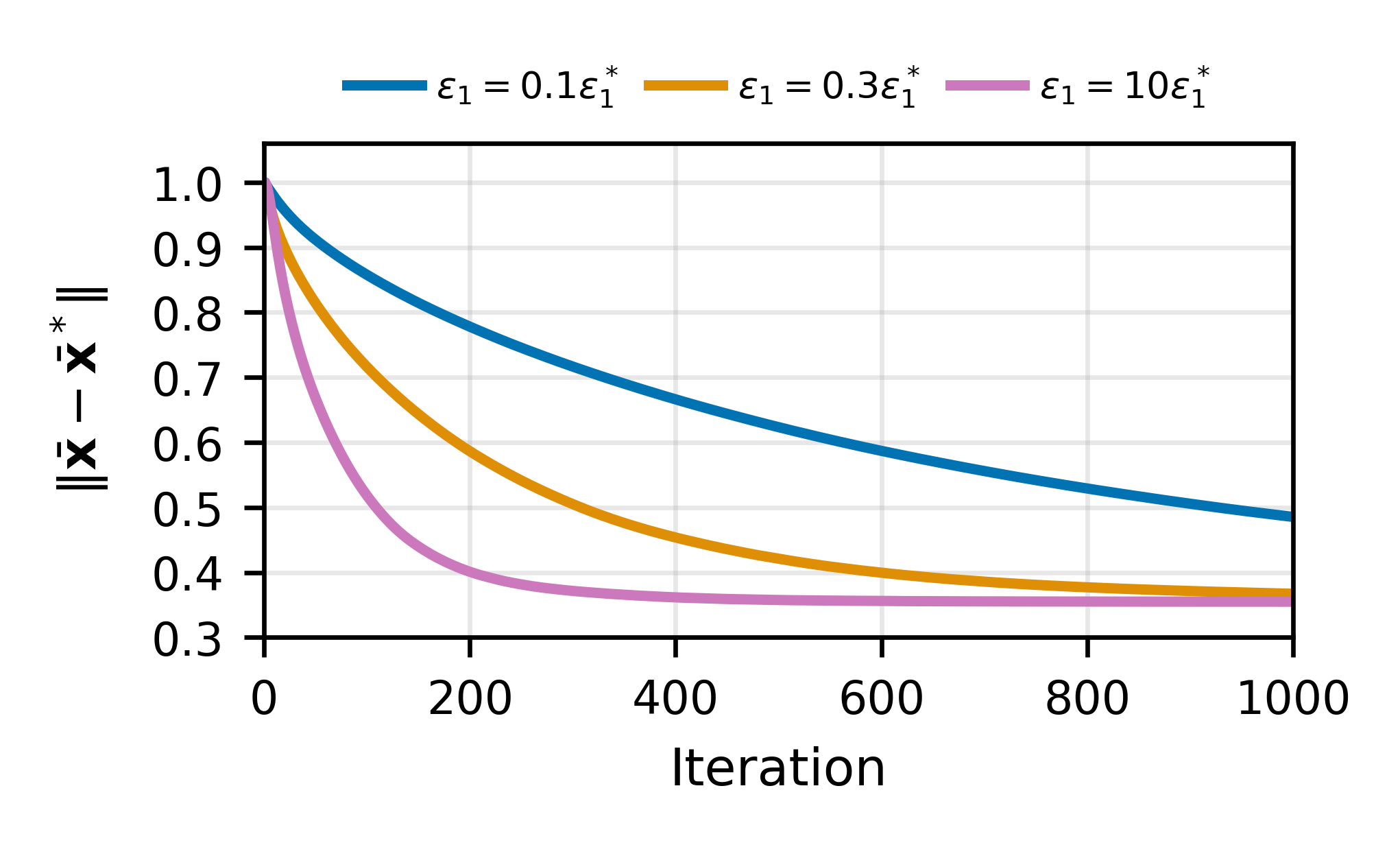}
\caption{Two-phase convergence trajectories for
three representative \(\epsilon_1\) values.}
\label{fig:rate_verification_a}
\end{figure}

Figure~\ref{fig:rate_verification_a} shows the
convergence trajectories for three representative
\(\epsilon_1\) values in log-linear scale.  Two
phases are visible, a rapid initial drop
(boundary layer) followed by a linear regime
(slow-phase exponential decay), as predicted by
singular perturbation theory.  Smaller
\(\epsilon_1\) yields a lower error floor but
slower convergence, while larger \(\epsilon_1\)
converges faster but plateaus at a higher level,
matching the rate--accuracy trade-off of
Theorem~\ref{thm:convergence}.

\begin{figure}[t]
\centering
\includegraphics[width=\columnwidth]{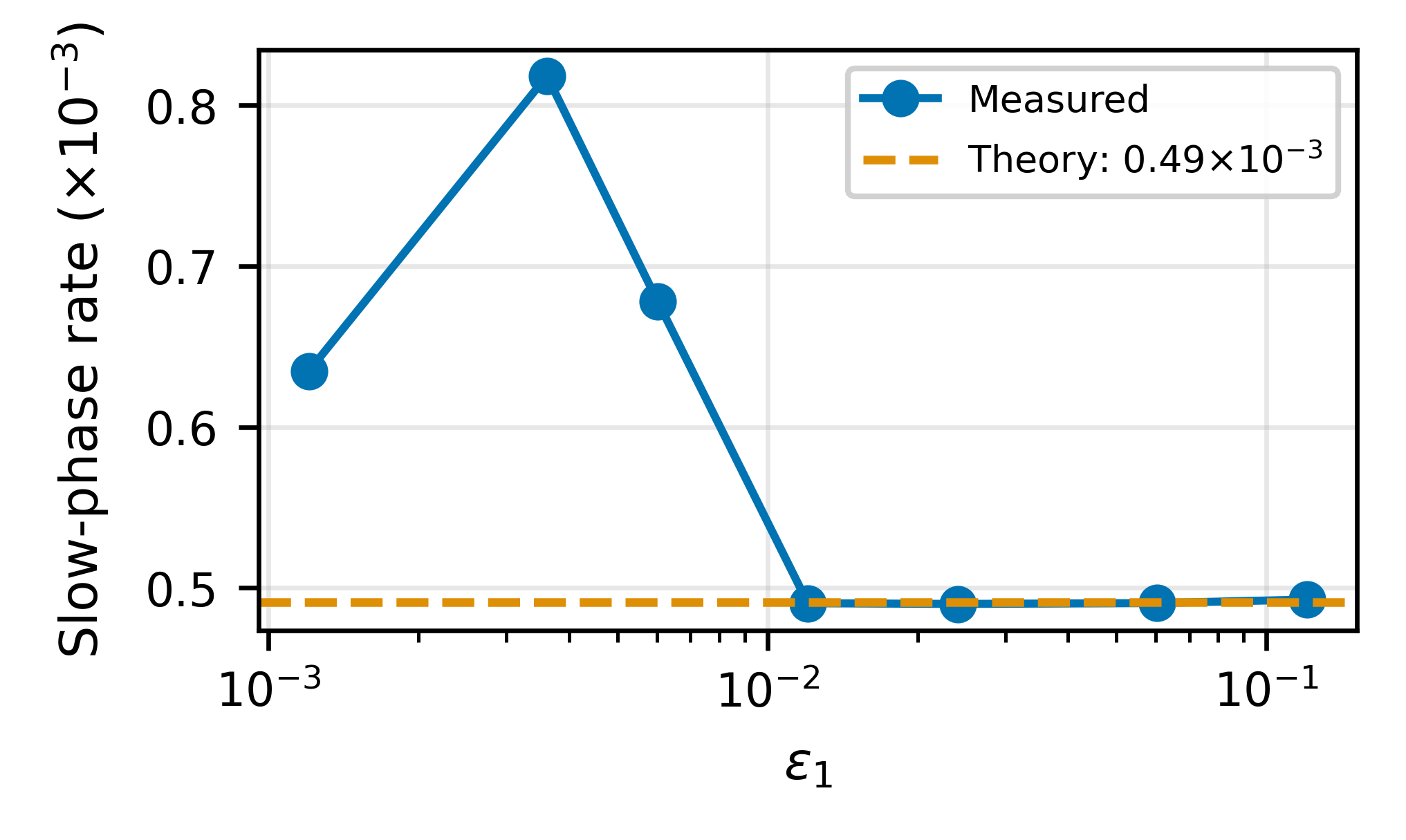}
\caption{\(\epsilon_1\)-independence of the
slow-phase convergence rate (\(\times 10^{-3}\)).}
\label{fig:rate_verification_b}
\end{figure}

Figure~\ref{fig:rate_verification_b} verifies
the \(\epsilon_1\)-independence of the slow-phase
rate.  Across the regime
\(\epsilon_1 \ge \epsilon_1^*\), the measured
rate remains nearly constant at approximately
\(0.49 \times 10^{-3}\) iterations\(^{-1}\),
with a coefficient of variation of only
\(0.2\%\), confirming that \(\zeta/(2k_1)\) is
governed by the slow subsystem and does not
depend on \(\epsilon_1\).  For
\(\epsilon_1 < \epsilon_1^*\), the rate is
slightly elevated
(\(0.6\)--\(0.8 \times 10^{-3}\)), reflecting
the stronger contribution of the boundary layer
transient.  At the three representative values
\(\epsilon_1/\epsilon_1^* = 0.1\), \(0.3\), and
\(10\), the measured rates are \(0.635\),
\(0.818\), and \(0.493 \times 10^{-3}\)
iterations\(^{-1}\), respectively, against the
theoretical \(0.49 \times 10^{-3}\).

\section{Conclusion}
\label{sec:conclusion}
\begingroup\sloppy\emergencystretch=3em
This paper proposed a bilevel distributed
PSO-neurodynamic method for solving the
distributed SCAD-regularized sparse
least-squares problem to global optimality.
A structure-preserving mollification of SCAD
restores differentiability without sacrificing
unbiasedness on large coefficients, so the
inner collaborative multi-agent neurodynamic
system converges using only local data and
neighbor-to-neighbor communication while
avoiding the second-order sufficient condition
invoked by existing distributed neurodynamic
methods, and the outer particle swarm search
mitigates the spurious local minima and
initialization sensitivity caused by the SCAD
nonconvexity.  A two-timescale singular
perturbation analysis yields a full-process
convergence characterization with explicit
exponential rates for the boundary-layer
transient and the slow phase, which in turn
enables the red-light--green-light mechanism
to freeze insensitive parameters and enforce
the singular-perturbation threshold.  Across
four sparse-recovery scenarios, the proposed
method achieves exact sparsity recovery with
faster convergence and lower relative error
than competing neurodynamic methods, and a
parameter sweep corroborates the predicted
rate--accuracy trade-off.
\par\endgroup

\end{document}